\documentclass[11pt]{amsart}
\usepackage{amsfonts,amssymb,amsmath}
\usepackage[latin1]{inputenc}
\usepackage[dvips]{epsfig}
\usepackage{graphicx}

\def\cal{\mathcal}
\newtheorem{theorem}{Theorem}[section]
\newtheorem{lemma}[theorem]{Lemma}
\newtheorem{assumption}[theorem]{Assumption}
\newtheorem{corollary}[theorem]{Corollary}
\newtheorem{proposition}[theorem]{Proposition}
\newtheorem{definition}[theorem]{Definition}
\newtheorem{remark}[theorem]{Remark}
\newcommand{\bq}{\begin{equation}}
\newcommand{\eq}{\end{equation}}
\newcommand{\bqa}{\begin{eqnarray}}
\newcommand{\eqa}{\end{eqnarray}}
\newcommand{\bd}{\begin{displaymath}}
\newcommand{\ed}{\end{displaymath}}
\def\beglem{\begin{lemma}\sl }
\def\endlem{\end{lemma}}
\def\begthm{\begin{theorem}\sl }
\def\endthm{\end{theorem}}
\def\begprop{\begin{proposition} \sl}
\def\endprop{\end{proposition}}
\def\begcor{\begin{corollary}\sl }
\def\endcor{\end{corollary}}
\def\begdef{\begin{definition}\sl}
\def\enddef{\end{definition}}
\def\begproof{ \noindent {\em Proof: \ }}
\def\endproof{\null\hfill {$\Box$}\bigskip}  
\def\beghyp{\begin{assumption} \sl}
\def\endhyp{\end{assumption}}
\def\begrem{\begin{remark}\rm}
\def\endrem{\null\hfill {$\Box$}\end{remark}}

 \def\C{{\cal C}}
 
 \def\E{{\cal E}}
 \def\F{{\cal F}}

 \def\Q{{\cal Q}}

 \def\T{{\cal T}}

 \def\ds{\displaystyle}
 \def\fref#1{(\ref{#1})}
 
 \newcommand{\RR}{\mathbb R}
 \newcommand{\NN}{\mathbb N}

 \def\null{\hbox{}}

 \def\half{1/2}
 
\def\medge{\meas(\sigma)}
 
 \def\({\left(}
 \def\){\right)}
 \def\[{\left[}
 \def\]{\right]}
 
 \def\2{{(2)}}

\def\Eint{\E_{int}}

\def\EextK{{\mathcal E}_{ext,K}}
\def\EintK{{\mathcal E}_{int,K}}

\def\meas{{\rm m}}
\def\dist{{\rm d}}

\def\beq{\begin{equation}}
\def\eeq{\end{equation}}
\def\ts{\tau_\sigma}

\def\signff{\bigskip\bigskip\hspace{80mm}
\vbox{{\sc Francis Filbet\par\vspace{3mm}
Universit\'e Lyon,  \par 
Universit\'e Lyon1, CNRS, \par  
UMR 5208 - Institut Camille Jordan, \par 
43, Boulevard du 11 Novembre 1918,\par 
F-69622 Villeurbanne cedex,
FRANCE\par\vspace{3mm}
e-mail:} filbet@math.univ-lyon1.fr }}

\begin{document}

\title[Asymptotically stable scheme for diffusive coagulation-fragmentation models]
{An asymptotically stable scheme for diffusive coagulation-fragmentation models}

\author{Francis Filbet}

\hyphenation{bounda-ry rea-so-na-ble be-ha-vior pro-per-ties
cha-rac-te-ris-tic}

\begin{abstract}
This paper is devoted to the analysis of a numerical scheme for the coagulation and fragmentation equation with diffusion in space. A finite volume scheme is developed, based on a conservative formulation of the space nonhomogeneous coagulation-fragmentation model, it is shown that the scheme preserves positivity, total volume and  global steady states. Finally, several numerical simulations are performed to investigate the long time behavior of the solution.
\end{abstract}

\maketitle




\section{Introduction}\label{sec:intro}
\setcounter{equation}{0}
Coagulation and fragmentation processes occur in the dynamics of cluster growth and describe the mechanisms by which clusters can coalesce to form larger clusters or fragment into smaller ones. Models of this type have found important applications in areas ranging over polymer kinetics \cite{ziff}, aerosols \cite{ref6}, cluster formation in astrophysics \cite{MFFK98} to the animal grouping in biology \cite{ref8}.

This paper is  devoted to the numerical simulation of the dynamics of the density function $f=f(t, x,y)\geq 0$ of particles (polymers, clusters) with position $x\in \Omega\subset \RR^N$ ($N\geq 1$), volume $y \in \RR^+:=(0,\infty)$, and time $t\geq 0$. The distribution function $f$ is  subject to coagulation and fragmentation phenomena with  respect to the volume variable $y$ and diffusion in space $x\in\Omega$ and is governed by the following equation
\begin{equation}
\label{full:model}
\frac{\partial f}{\partial t} \,-\, d(y)\, \Delta_x f \,= \Q(f), 
\end{equation}
where $\Q(f)$ is the coagulation-fragmentation operator acting on the volume variable $y\in \RR^+$ and these particles diffuse in the environment $\Omega$, which is assumed to be a smooth bounded domain with normalized volume $|\Omega|=1$. Moreover, the model (\ref{full:model}) is supplemented with an initial datum 
\begin{equation}
\label{hypf0}
f(0,x,y) := f^{\rm in}(x,y),
\end{equation}
and in the following analysis we will consider homogeneous Neumann boundary conditions on $\partial \Omega$
\begin{equation}
\label{hypBC}
\nabla_x f(t,x,y)\,\cdot \, \nu(x) \,=\, 0,   \quad (t,x,y)\in \RR^+ \times \partial \Omega \times \RR^+,
\end{equation}
with $\nu$ the  outward unit normal to $\Omega$. We also assume the diffusion coefficient $d(y)$ to be non degenerate in the sense that there exists $d^*, d_* \in \RR^+$ such that
\begin{equation}
d^* \geq d(y) \geq d_* >0, \quad y \in \RR^+.
\label{hypd}
\end{equation}    
Let us first describe  precisely the coagulation and fragmentation process.  In the simplest coagulation-fragmentation models the clusters are usually assumed to be fully identified by their size (or volume, or number of particles). The coagulation-fragmentation models we consider in this paper describe the time evolution of the cluster size distribution as the system of clusters undergoes binary coagulation and binary fragmentation events. More precisely, denoting by $C_y$ the clusters of size $y$ with $y \in\RR^+$, the basic reactions taken into account herein are 
\begin{equation}
C_{y} + C_{y'}\;\;\; {{a(y,y')}\atop{\longrightarrow}}\;\;\; C_{y+y'},\;\;\;\mbox{(binary coagulation)} 
\label{particle1}
\end{equation}
and
\begin{equation}
C_{y} \;\;\; {{b(y-y',y')}\atop{\longrightarrow}} \;\;\; C_{y-y'} + C_{y'},\;\;\;\mbox{(binary fragmentation)},
\label{particle2}
\end{equation}
where $a$ and $b$ denote the coagulation and fragmentation rates respectively, and are assumed to depend only on the size of the clusters involved in these reactions.

Thus, the coagulation-fragmentation operator is defined by
\begin{equation}
\label{Qf}
\Q(f) \,=\, \Q_\texttt{c}(f) - \Q_\texttt{f}(f),
\end{equation}
where the coagulation operator is 
\begin{equation}
\Q_\texttt{c}(f)(x,y)  \,= \, \frac{1}{2}\int_0^y a(y', y-y')\, f(x,y') \,f(x,y-y')\, dy'\, \, - \,\, \int_0^\infty a(y,y')\, f(x,y)\,f(x,y')\, dy' 
\label{Qcf} 
\end{equation}
whereas fragmentation mechanism by which a single particle splits into two pieces is given by
\begin{equation}
\Q_\texttt{f}(f)(x,y) \, = \, \frac{1}{2} \int_0^y b(y', y-y')\, dy'\, f(x,y)\, \, - \,\,\int_0^\infty b(y,y')\, f(x,y+y')\, dy'. 
\label{Qff}
\end{equation}
The coagulation coefficient, $a=a(y,y')$, characterizes the rate at which the coalescence of two particles with respective volumes $y$ and $y'$ produces a particle of volume $y+y'$, whereas the fragmentation coefficient, $b=b(y,y')$, represents the rate at which the fragmentation of one particle with volume $y+y'$ produces two particles of volume $y$ and $y'$. Both coefficients $a$ and $b$ are nonnegative  functions and 
\begin{equation}
\label{hypa}
\left\{
\begin{array}{l}
\displaystyle{a(y,y')=a(y',y), \,\, b(y,y')=b(y',y),} 
\\ 
\,
\\
\displaystyle{a, \,\, b \,\in L^\infty_{loc}\left(\RR^+\times \RR^+ \right).}
\end{array}\right.
\end{equation}

For symmetric kernels, we observe that during the microscopic coagulation and fragmentation processes, as depicted in equations (\ref{particle1})-(\ref{particle2}), the number of particles varies with time while the total volume of particles is conserved.

In terms of $f$, the total number of particles and the total volume of particles
at time $t\geq 0$ are respectively  given by
$$
M_0(t,x) := \int_{\RR^+} f(t,x,y) \,dy,\quad M_1(t,x) := \int_{\RR^+} y\, f(t,x,y) \,dy. 
$$
Besides existence and uniqueness results, few is known on the qualitative behavior of solutions of coagulation-fragmentation equations except when coagulation and fragmentation coefficients are linked by the so called ``detailed balanced condition'' : there exists a nonnegative function $M\in L^1 (\Omega\times \RR^+,(1+y)dxdy)$ such that
\begin{equation}
\label{dbc}
a(y,y') \, M(y)\, M(y') \,\,=\,\, b(y,y') \, M(y+y'), \quad (x,y,y')\in\Omega \times \RR^+\times \RR^+.
\end{equation}
This condition implies that $M$ is a stationary solution to the coagulation-fragmentation equation 
$$
d(y) \,\Delta_x f \,+\, \Q(f) \,=\, 0 .
$$ 
We further assume that, for each $R\in\RR^+$, the equilibrium $M$ satisfies the positivity condition
\begin{equation}
\label{hypM}
\inf_{y\in[0,R]} M(y) \,>\,0.
\end{equation} 

An additional and interesting feature of the detailed balanced condition is that there exists an entropy $H$ given by
$$
H(f|M) \,=\, \int_{\Omega\times \RR^+} \left[f(x,y) \left( \ln\left(\frac{f(x,y)}{M(y)}\right)\,-\, 1 \right)\,+\, M(y)\right]\, dy\,dx, 
$$
which satisfies the following $H$ theorem
$$
\frac{dH}{dt} \,=\, -\frac{1}{2}\, D(f)  \,-\, \int_{\Omega\times \RR^+} \frac{\left( \nabla_x f\right)^2}{f}(t,x,y)\,dy\,dx\,\leq\, 0, 
$$  
where the entropy dissipation of the coagulation-fragmentation operator is defined by 
$$
D(f) \,=\,\int_{\Omega}\int_{\RR^+\times \RR^+} \left(a\,f\,f'\, -\,  b \,f'' \right)\,\left( \ln(a\,f\,f') - \ln(b\,f'')\right)\,dy\,dy'\,dx,  
$$
with the shorthand notation $f:=f(t,x,y)$, $f':=f(t,x,y')$ and $f'':=f(t,x,y+y')$.  Since $D(f)$ only vanishes when $f$ is an equilibrium, it naturally means that
$$
f(t) \longrightarrow M, \quad \textrm{when }\quad t \rightarrow \infty. 
$$ 

The main purpose of this work is to present a numerical scheme to solve (\ref{full:model}) built upon a finite volume discretization 
with respect to the space variable $x\in\Omega$ and  volume   variable $y\in \RR^+$. The analysis of the so-obtained scheme would allow to prove the convergence of the discretized particle density towards a solution to the continuous problem on a fixed time interval $[0,T]$ ($T>0$), see for instance \cite{bourg:fil}, but here our  aim is different. It will consist in the study of the long time behavior of the numerical solution  on a fixed mesh in space and volume variables. Indeed, often in applications we are interested in steady states or in the long time behavior of the solution, it is then important to design a numerical scheme, which has good stability properties uniformly in time and remains consistent with respect to the exact solution for long time. Thus, the scheme proposed in this paper is designed such that it preserves qualitative properties of the exact solution as  steady states of the coagulation-fragmentation operators  (\ref{Qcf})-(\ref{Qff}). Moreover, an estimate of the entropy dissipation due to an appropriate finite volume approximation with respect to space and volume variables  is given and the scheme is shown to give a consistent  approximation of the continuous problem in the  long time asymptotic limit $t\rightarrow \infty$.

Before describing more precisely our results, let us recall that the
coagulation and fragmentation equations (\ref{full:model}) with (\ref{Qcf})-(\ref{Qff}) have been the object of several studies recently. 

On the one hand, among the various approaches for the approximation of coagulation and fragmentation models, we may distinguish between deterministic and Monte Carlo methods. We refer for instance to \cite{EEE94,FLxx,Le00} for deterministic methods, \cite{Ba99, guias} for stochastic methods, and the references therein.  Concerning the convergence analysis of numerical methods for coagulation and fragmentation models, we refer to \cite{LW04} for a rigorous study of quasi Monte-Carlo methods. For deterministic approximations, the situation is different since the relationship between  discrete and  continuous models has been considered by some authors, see the survey paper \cite{Dr72} and \cite{aizen}. 
A rigorous setting for the formal analysis performed in \cite{aizen} under general assumptions on the coagulation and fragmentation coefficients has been given in \cite{LM02}. Then, similar techniques are used to prove convergence of discrete schemes to the exact solution in \cite{bourg:fil}.

However, few is known concerning the stability and the analysis of the numerical solution in the long time asymptotic limit.  For the continuous model, in reference \cite{LM01}, the authors are able to apply the techniques of weak compactness, as in the paper by R.J. DiPerna and P.-L. Lions \cite{dpl} about the Boltzmann equation, to prove the weak stability of weak solutions. In the case when an $H$-theorem holds, they can obtain some partial information about the large-time behavior of solutions. More recently, J. A. Carrillo {\it et al.} prove  exponential decay of the solution  towards equilibrium for the inhomogeneous Aizenman-Bak model \cite{CDF} using entropy dissipation methods \cite{DV}.  Our analysis for discrete models will be inspired by these works. 


We now briefly outline the contents of the paper. In the next section,
we introduce the numerical approximation of (\ref{full:model}) and state the stability result which we prove in Section \ref{apriori::sec}.  In the final section (Section~\ref{test:num}), some numerical simulations are performed with the numerical scheme presented in Section~\ref{sec2} and the long time behavior of the solution is investigated.

\section{Numerical scheme and main results} \label{sec2}
In order to compute an approximation of this model using a finite volume method in space variable $x\in\Omega$ and volume variable $y\in\RR^+$, we observe that the coagulation-fragmentation operator can be written in a conservative form. Indeed, writing equation \fref{full:model} in a ``conservative'' form, as proposed in \cite{MFFK98, TIN96}, enables to describe precisely the time evolution of the total volume. Also, this formulation is particularly well adapted to a finite volume discretization 
which, in turn, is expected to give a precise account of volume conservation. Precisely, the coagulation and fragmentation terms can be written in divergence form:
\begin{eqnarray*}\left\{
\begin{array}{l}
\displaystyle{y\,\Q_\texttt{c}(f)(x,y)  =  - \frac{\partial \C(f)}{\partial y}(x,y), }
\\
\, 
\\
\displaystyle{y\,\Q_\texttt{f}(f)(x,y) =  - \frac{\partial \F(f)}{\partial y}(x,y),}
\end{array}\right.
\end{eqnarray*}
where the operator $\C(f)$ is given by 
\begin{eqnarray}
\label{coag:op}
\C(f)(x,y) &:=& \int_0^y \int_{y-u}^\infty u\ a(u,v)\ f(x,u)\ f(x,v)\ dvdu\,,
\;\;\; (x,y)\in\Omega\times\RR^+, 
\end{eqnarray}
and $\F(f)$ is
\begin{eqnarray}
\label{frag:op}
\F(f)(x,y) &:=& \int_0^y \int_{y-u}^\infty u\ b(u,v)\ f(x,u+v)\ dvdu\,,
\;\;\; (x,y)\in\Omega\times\RR^+. 
\end{eqnarray}
Finally, the coagulation-fragmentation equation  reads
\begin{eqnarray}
\label{evol:eq}
\left\{
\begin{array}{l}
\displaystyle{y\ \frac{\partial f}{\partial t} \,-\, d(y)\,y\,\Delta_x f \,=\, - \frac{\partial\, \C(f)}{\partial y} \,\,+\,\, \frac{\partial\, \F(f)}{\partial y}\,,} 
\\
\,
\\
f(0,x,y) = f^\texttt{in}(x,y),  \;\;\; (x,y)\in\Omega,\times\RR^+
\end{array}
\right.
\end{eqnarray}
and we assume that the initial datum $f^\texttt{in}$ is a nonnegative function which satisfies: 
\begin{equation}
\label{hypL1}
f^\texttt{in} \in L^{1}(\Omega\times\RR^+)\cap L^1(\Omega\times\RR^+,y dxdy).
\end{equation}
Here and below, the notation $L^1(\Omega\times\RR^+,y dxdy)$ stands for the space of the Lebesgue measurable real-valued functions on $\Omega\times \RR^+$ which are integrable with respect to the measure $y dxdy$.

When designing the volume discretization 
of the coagulation and fragmentation terms, it is necessary to truncate the infinite integrals in formulae \fref{coag:op}-\fref{frag:op}. But this means restricting the domain of action of kernels $a$ and $b$ to a bounded set of volumes $y$, that is, preventing coagulation to occur among particles with volume exceeding a fixed value. The discretization 
we propose is based on the conservative truncation method for the coagulation and fragmentation terms. Given a positive real $R$, let $(x,y)\in \Omega\times(0,R)$
\begin{eqnarray}
\nonumber
\C^R(f)(x,y)  &:=& \int_0^y \int_{x-u}^{R-u} u\ a(u,v)\ f(x,u)\ f(x,v)\ dvdu\,,
\\
&=& \int_0^y \int_{y}^{R} u\ a(u,w-u)\ f(x,u)\ f(x,w-u)\ dw\,du\,.
\label{trunc:coag}
\end{eqnarray}  
In that case, $\C^R(f)(x,0)\,=\,\C^R(f)(x,R) \,=\, 0$ foreach $x\in\Omega$ so that the total volume of the solution is now nonincreasing with respect to time. 

As regards the fragmentation term, the truncation is also a conservative truncation on the fragmentation term. Using the same idea, we introduce for  $(x,y)\in \Omega\times(0,R)$
\begin{eqnarray}
\nonumber
\F^R(f)(x,y) &:=& \int_0^y \int_{y-u}^{R-u} u \,b(u,v) \,f(x,u+v) \,dv\, du. 
\\
&=& \int_0^y \int_{y}^{R} u \,b(u,w-u) \,f(x,w) \,dv\, du. 
\label{frc}
\end{eqnarray}
Then, the conservative coagulation-fragmentation operator satisfies exactly the conservation of total volume, so that the following equation is indeed a truncated conservative coagulation and fragmentation equation:
\begin{eqnarray}
\left\{
\begin{array}{l}
\displaystyle{y \frac{\partial f_R}{\partial t} \,-\, d(y)\,y\,\Delta_x f_R \,=\, -\,\frac{\partial \C^R(f_R)}{\partial y}(x,y) \,\,+\,\, \frac{\partial \F^R(f_R)}{\partial y}(x,y)\,} 
\label{evian} 
\\
\,
\\
f(0,x,y) \,=\, f^\texttt{in}(x,y),  \quad (x,y)\in\Omega \times (0,R),  
\end{array}
\right.
\end{eqnarray}
since 
$$
\frac{d}{dt}\int_0^R\int_{\Omega} y \,f_R(t,x,y)\, dx\,dy \,=\,  0.
$$
Under the detailed balance condition (\ref{dbc}), model (\ref{evian}) has also a steady state $M_R$, only depending on $y$ and such that 
$$
M_R(y) \,=\, M(y), \quad y \in [0,R].
$$
and 
$$
\frac{1}{|\Omega|}\,\int_0^R\int_{\Omega} f^\texttt{in}(x,y)\,dx\,dy \,=\, \int_0^R M_R(y)\,dy.
$$
Convergence for large values of $R$ has been thoroughly studied in the recent past. We briefly mention some results for the coagulation equation (that is, with $b=0$). These results adapt easily to the coagulation-fragmentation equation but under different assumptions on the kernels. 
Let us mention that when $a(y,y')/(y\ y')\to 0$ as $y+y'\to +\infty$, convergence as $R\to +\infty$ of the solutions to (\ref{evian}) toward a solution of (\ref{coag:op})-(\ref{evol:eq}) can be proven by using the approach  developed in \cite{LM02}. Thus, since the convergence of solutions to \fref{evian} towards solutions of (\ref{coag:op})-(\ref{evol:eq}) is well established in rather general situations, this paper will only focus on the convergence of a sequence built on a numerical scheme towards a solution to the equation \fref{evian} when the truncature $R$ is fixed. 
In the remainder of the paper, for the sake of clarity, we drop the subscript $R$ and write $f$ instead of $f_R$ for a solution of equation \fref{evian}. Parameter $R$ being fixed, this should raise no confusion.

Now, we turn to the discretization 
of equation \fref{evian}. Having reduced the computation to a bounded interval, the second step
is to introduce the space and  volume discretizations. 
To this end, we set $\Delta y\in (0,1)$ and $N_y$  a large integer, and denote by 
$(y_{i-1/2})_{i\in\{0,\ldots,N_y\}}$ a mesh of $(0,R)$, where $y_{i-\half} \,\,=\,\, i\,\Delta y$, and $\Lambda_i=[y_{i-\half},y_{i+\half})$ for $i\ge 0$. 

Concerning the space variable, we choose an admissible mesh of $\Omega$ given by a family $\T$ of control volumes (open and convex polygons in 2-D), a family $\E$ of edges and a  family of points $(x_K)_{K\in\T}$ which satisfy Definition 5.1 in \cite{egh}. It implies that the straight line between two neighboring centers of cells $(x_K,x_L)$ is orthogonal to the edge $\sigma =K|L$. In the set of edges $\E$, we distinguish the interior edges $\sigma\in \E_{int}$ and the boundary edges $\sigma\in \E_{ext}$. For a control volume $K\in\T$, we denote by $\E_K$ the set of its edges, $\EintK$ the set of its interior edges, $\EextK$ the set of edges of $K$ included in $\Gamma=\partial\Omega$. 

In the sequel, we denote by $\dist$ the distance in $\RR^N$, $\meas$ the measure in $\RR^N$. We assume that the family of mesh considered satisfies the following regularity constraint :  there exists $\xi >0$  such that 
\begin{equation}
\label{mesh00}
\dist(x_K,\sigma) \,\,\geq\,\,  \xi\, \dist(x_K,x_L),\quad  \mbox{ for } K\in \T, \mbox{ for } \sigma \,\in\EintK , \ \sigma= K|L.
\end{equation} 
The size of the mesh is defined by
\beq\label{defh}
\delta\,=\,\ds\max_{K\in\T}\left(\rm{diam}(K)\right).
\eeq
For all $\sigma\in\E$, we define the transmissibility coefficient:
$$
\ts=\left\{\begin{array}{ll}
\ds\frac{\medge}{\dist(x_K,x_L)},&\mbox{ for }\sigma\in \E_{int},\ \sigma=K|L,\\[5.mm]
\ds\frac{\medge}{\dist(x_K,\sigma)},&\mbox{ for }\sigma\in\EextK. \\
\end{array}
\right.
$$
Next, $X(\mathcal{T})$ will be the set of functions from $\Omega$ to $\RR$ which are constant over each control volume $K\in\mathcal{T}$.

We define the approximation $f^{h}(0)$ of the initial datum $f^{\rm in}$ as usual by
\begin{equation}
\label{discrinit}
f^{h}(0) \,=\, \sum^{N_y}_{i=0}\sum_{K\in\mathcal{T}} f_{K,i}^{\rm in} \,\mathbf{1}_{K\times\Lambda_{i}},
\end{equation}
with 
$$
f_{K,i}^{\rm in} = \frac{1}{{\rm m}(K)\,\Delta y}\ \int_{\Lambda_{i}}\int_{K} f^{\rm in}(x,y)\ dx\,dy\,, 
$$
where $\mathbf{1}_E$ denotes the characteristic function of the subset $E$ of $\Omega$ and converges strongly to $f^0$ in $L^1(\Omega\times (0,R))$ as $h=(\Delta y,\delta)$ goes to $0$.

Let us now introduce the numerical scheme itself. For each integer $i\in\{0,\cdots,N_y\}$ and each $K\in\mathcal{T}$, we define the approximation of $f(t, x,y)$ for $t\in\RR^+$ and $(x,y)\in K \times \Lambda_i$ as $f_{K,i}(t)$. The sequence $(f_{K,i})_{K,i}$ is defined by the following discretization 
of the coagulation-fragmentation equation\,: for  $K\in\mathcal{T}, \ i\in\{0,\hdots,N_y\}$, we  solve the ordinary differential system
\begin{equation}
\left\{
\begin{array}{l}
\displaystyle{{\rm m}(K) \, \frac{df_{K,i}}{dt}  \,-\, d(y_{i-\half})\sum_{\sigma\in \mathcal{E}_K}\tau_{\sigma} \,D_{K,\sigma}f_{K,i} \, =\,  {\rm m}(K) \,\Q_{K,i},} 
\\
\,
\\
\displaystyle{f_{K,i}(0) := f_{K,i}^{\rm in}.}
\end{array}\right.
\label{discr}
\end{equation}
We have set
\begin{equation}
\label{QKi}
 \Q_{K,i} \,=\,
-\frac{\C_{K,i+1/2} - \C_{K,i-1/2} }{y_{i-\half}\,\Delta y} + \frac{\F_{K,i+1/2} - \F_{K,i-1/2} }{y_{i-\half}\,\Delta y} 
\end{equation}
where the flux $\C_{K,i+1/2}$ in (\ref{QKi}) represents an approximation of the coagulation operator (\ref{trunc:coag}), $\F_{K,i+1/2}$ is the approximation of the fragmentation part (\ref{frc}) and both are defined by
\begin{eqnarray}
 \C_{K,i+1/2} &=& \sum^{i}_{j=0}\sum^{N_y-1}_{l=i+1} \Delta y^2 \,y_{j-\half} \,a_{j,l-j}\, f_{K,j} \, f_{K,l-j}, 
\label{discrC} 
\end{eqnarray}
where $a_{i,j}\,:=\,a(y_i,y_j)$ with $y_i=(i+1/2)\,\Delta y$ and
\begin{eqnarray}
 \F_{K,i+1/2} &=& \sum^{i}_{j=0}\sum^{N_y-1}_{l=i+1} \Delta y^2 \,y_{j-\half}\,b_{j,l-j}\,f_{K,l}, 
\label{discrF} 
\end{eqnarray}
where $b_{i,j}\,:=\,b(y_i,y_j)$ whereas the fluxes at the boundary are
\begin{eqnarray}	
\C_{K,-1/2} = \F_{K,-1/2} =\C_{K,N_y+1/2}  =  \F_{K,N_y+1/2} = 0,\quad K\in\mathcal{T}.
\label{discrboundary}
\end{eqnarray}
Concerning the approximation of the diffusion in space, we set $\nu_{K,\sigma}$ the unit normal to $\sigma$ outward from $K$ and define an approximation of $\nabla_x f \cdot \nu_{K_\sigma}$ on $\sigma$ by
\begin{equation}
\label{schem3}
D_{K,\sigma}f_{K,i} =  \left\{
\begin{array}{ll}
\displaystyle{f_{L,i} - f_{K,i}, \textrm{ if }   \sigma=K|L\in \mathcal{E}_{int,K},}
\\
\displaystyle{0,\textrm{ if }   \sigma\in \mathcal{E}_{ext,K},}
\end{array}\right.
\end{equation}
for all $K \in \mathcal{T}$. This discretization obviously relies on  a simple finite volume approach for the space and volume variables (see, {\em e.g.}  \cite{bourg:fil}).

The following function $f^h$ defined on $\RR^+\times\Omega\times [0, R]$ will be useful in the sequel.
\begin{equation}
f^h(t, x,y) = \sum_{K\in\mathcal{T}}\sum^{N_y}_{i=0} f_{K,i}(t)\,\mathbf{1}_{K\times\Lambda_{i}}.
\label{fhdef}	
\end{equation}

We may now state our main result.

\begthm
Assume that the coagulation and fragmentation kernels satisfy (\ref{hypd}), (\ref{hypa}) and (\ref{hypM}) and $f^\texttt{in}$ satisfies (\ref{hypL1}). We consider a uniform volume mesh in $y$ and require the mesh $\T$ in space  to satisfy condition (\ref{mesh00}), (\ref{defh}).

Then,  there exists a unique solution  $f^h$ to the finite volume scheme (\ref{discr})-(\ref{schem3}), which satisfies
\begin{itemize}
\item[(i)] nonnegativity  
$$
\displaystyle{f_{K,i}(t) \geq 0, \quad  (K,i) \in \mathcal{T} \times \{0,\ldots,N_y\}, \quad  \, t \in \RR^+,}
$$
\item[(ii)] conservation of total volume
$$
\displaystyle{\frac{d}{dt} \left[\sum_{K,i} {\rm m}(K) \, \Delta y\, y_{i-1/2}\,f_{K,i}(t)\right] = 0, \quad \forall \, t \in \RR^+,}
$$
\item[(iii)] entropy dissipation
$$
\displaystyle{\frac{dH(f^h|M) }{dt} \,+\,  \Delta y\,\sum_{K,i} \,\sum_{\sigma\in \mathcal{E}_K}d(y_{i-\half}) \tau_{\sigma} \,D_{K,\sigma}f_{K,i} \ln\left(\frac{f_{K,i}}{M_i}\right) \,=\, -\,D(f^h),}
$$
\end{itemize}
where 
$$
H(f^h|M) \,=\, \sum_{K,i} {\rm m}(K)\, \Delta y\,\left[ f_{K,i}\,\left(\ln\left(\frac{f_{K,i}}{M_i}\right)\,-\, 1  \right)  + M_i\right] 
$$
and 
$$
D(f^h) \,=\, \frac{1}{2}\sum_{K,i}\sum^{i}_{l=0} \Delta y^2 \,{\rm m}(K)\,\left( a_{l,i-l}\, f_{K,l} \, f_{K,i-l} \,-\, b_{l,i-l}\,f_{K,i} \right)\left[ \ln\left(\frac{a_{l,i-l}\, f_{K,l} \, f_{K,i-l}}{b_{l,i-l}\,f_{K,i}}\right) \right].
$$
\begin{itemize}
\item[(iv)] Moreover, the finite volume scheme is asymptotic preserving {\it i.e.} 
\end{itemize}
$$
f^h(t) \longrightarrow M^h, \quad \textrm{ as } t \longrightarrow \infty.
$$
where
$$
M^h(x,y) \,=\, M_i \,:=\,\exp(-\alpha \, y_{i}), \quad\forall (x,y) \in K\times [y_{i-1/2},y_{i+1/2})  
$$ 
and $\alpha\in\RR^+$ is such that
$$
\sum_{K\in\T} \sum_{i=0}^{N_y} {\rm m}(K) y_{i-1/2} M_i \,=\, \sum_{K\in\T} \sum_{i=0}^{N_y} {\rm m}(K) y_{i-1/2} f^\texttt{in}_{K,i}.
$$
\label{thm01}
\endthm

\section{A priori estimates}\label{apriori::sec}
As we mentioned in the introduction, our goal is not to prove that the sequence of functions $(f^h)_h$ converges in some sense to a function $f$ as $h=(\Delta y,\delta)$ goes to $0$, but to study the long time behavior of the numerical solution and to prove its uniform stability with respect to time. First, we prove that the solution $f^h$ to the scheme (\ref{discr})-(\ref{discrboundary}) enjoys properties similar to those of function $f$ given by (\ref{evian}) which we gather in Proposition~\ref{prop:weak} below.

On the one hand, we first re-write the discrete coagulation-fragmentation obtained form the conservative form (\ref{evian}) in a new form which is consistent with the classical formulation (\ref{Qf})-(\ref{Qff}). On the other hand, we find a discrete version of the weak formulation from which we prove {\it a priori} estimates.
  
\begin{proposition}
\label{prop:reform}
Assume the approximation of the coagulation-fragmentation model (\ref{trunc:coag})-(\ref{evian}) is given  by the finite volume scheme  (\ref{discr})-(\ref{schem3}). Then, the discrete operator satisfies 
\begin{eqnarray*}
 \Q_{K,i} \,=\,
-\frac{\C_{K,i+1/2} - \C_{K,i-1/2} }{y_{i-\half}\,\Delta y} + \frac{\F_{K,i+1/2} - \F_{K,i-1/2} }{y_{i-\half}\,\Delta y} 
\end{eqnarray*}
with
\begin{eqnarray}
 \label{discrQ} 
\Q_{K,i} &=& \frac{1}{2}\sum^{i}_{j=0} \Delta y \,\left(a_{j,i-j}\,f_{K,j}\,f_{K,i-j}  \,-\, b_{j,i-j}\,f_{K,i}\right)
\\
&-& \sum^{N_y-1}_{j=i} \Delta y \,\left(a_{i,j-i}\,f_{K,i}\,f_{K,j-i}  \,-\, b_{i,j-i}\,f_{K,j}\right)
\nonumber
\end{eqnarray}
\end{proposition}
\begproof
Starting from the definition of the fluxes (\ref{discrC}) and (\ref{discrF}), we set 
$$
\Q_{K,i} \,=\,
-\frac{\C_{K,i+1/2} - \C_{K,i-1/2} }{y_{i-\half}\,\Delta y} + \frac{\F_{K,i+1/2} - \F_{K,i-1/2} }{y_{i-\half}\,\Delta y} 
$$
and by construction
\begin{eqnarray*}
y_{i-\half}\,\Q_{K,i} &=& -\sum^{i}_{j=0}\sum^{N_y-1}_{l=i+1} \Delta y \,y_{j-\half} \,\left( a_{j,l-j}\, f_{K,j} \, f_{K,l-j} \,-\, b_{j,l-j}\,f_{K,l} \right) 
\\
&+& \sum^{i-1}_{j=0}\sum^{N_y-1}_{l=i} \Delta y \,y_{j-\half} \,\left( a_{j,l-j}\, f_{K,j} \, f_{K,l-j} \,-\, b_{j,l-j}\,f_{K,l} \right), 
\end{eqnarray*}
which can be easily simplified into
\begin{eqnarray*}
y_{i-\half}\,\Q_{K,i} &=&  -\sum^{N_y-1}_{l=i+1} \Delta y \,y_{i-\half} \,\left( a_{i,l-i}\, f_{K,i} \, f_{K,l-i} \,-\, b_{i,l-i}\,f_{K,l} \right) 
\\
&+& \sum^{i-1}_{l=0} \Delta y \,y_{l-\half} \,\left( a_{l,i-l}\, f_{K,l} \, f_{K,i-l} \,-\, b_{l,i-l}\,f_{K,i} \right).
\end{eqnarray*}
Then, adding and subtracting the term $y_{i-\half}\left( a_{i,0}\, f_{K,0} \, f_{K,0} \,-\, b_{i,0}\,f_{K,i}\right)$, we get
\begin{eqnarray*}
y_{i-\half}\,\Q_{K,i}&=& \sum^{i}_{l=0} \Delta y \,y_{l-\half} \,\left( a_{l,i-l}\, f_{K,l} \, f_{K,i-l} \,-\, b_{l,i-l}\,f_{K,i} \right)
\\
&-&\sum^{N_y-1}_{l=i} \Delta y \,y_{i-\half} \,\left( a_{i,l-i}\, f_{K,i} \, f_{K,l-i} \,-\, b_{i,l-i}\,f_{K,l} \right). 
\end{eqnarray*}
Hence, it only remains  to treat the first term of the right hand side:  we compute for a symmetric sequence $g_{i,j}$ such that $g_{i,j}=g_{j,i}$ for each $(i,j)\in\{0,\cdots,N_y\}^2$
\begin{eqnarray*}
\sum_{l=0}^i l \,g_{l,i-l} &=& \sum_{l=0}^i i \,g_{l,i-l} \,-\,   \sum_{l=0}^i (i-l) \,g_{l,i-l}
\\
 &=& \sum_{l=0}^i i \,g_{l,i-l} \,-\,   \sum_{l=0}^i l \,g_{i-l,l}.
\end{eqnarray*}
Thus, using the symmetry $g_{i-l,l}=g_{l,i-l}$ we finally get
$$
\frac{1}{2}\sum_{l=0}^i i \,g_{l,i-l} \,=\, \sum_{l=0}^i l \,g_{l,i-l}.
$$
Therefore, with $y_{i-\half}=i\,\Delta y$ and taking successively $g_{i,j} =  a_{i,j}\, f_{K,i} \, f_{K,j}$ and $g_{i,j} =  b_{i,j}$, it yields the result
\begin{eqnarray*}
\Q_{K,i} &=& \frac{1}{2}\sum^{i}_{l=0} \Delta y \,\left( a_{l,i-l}\, f_{K,l} \, f_{K,i-l} \,-\, b_{l,i-l}\,f_{K,i} \right)
\\
&-&\sum^{N_y-1}_{l=i} \Delta y \,\left( a_{i,l-i}\, f_{K,i} \, f_{K,l-i} \,-\, b_{i,l-i}\,f_{K,l} \right). 
\end{eqnarray*}
\endproof

Next, we prove the following estimates

\begin{proposition}
\label{prop:weak}
Assume the approximation of the coagulation-fragmentation operator  (\ref{trunc:coag})-(\ref{frc}) is given  by the finite volume scheme  (\ref{discrC})-(\ref{discrF}). Then, for each $K\in\T$ the discrete operator $(\Q_{K,i})_{0\leq i \leq N_y}$ satisfies a discrete weak formulation : for any sequence $(\varphi_{i})_{0\leq i \leq N_y}$ and for each $K\in\T$
\begin{equation}
\label{weak:d}
\sum_{i=0}^{N_y-1} \Delta y\,\Q_{K,i}\,\varphi_{i} \,=\, -\frac{1}{2}\sum_{i=0}^{N_y-1}\sum^{i}_{l=0} \Delta y^2 \,\left( a_{l,i-l}\, f_{K,l} \, f_{K,i-l} \,-\, b_{l,i-l}\,f_{K,i} \right)\left(\varphi_{l} \,+\,\varphi_{i-l} \,-\, \varphi_{i} \right)
\end{equation}
\end{proposition}
\begproof
We multiply (\ref{discrQ}) by $\Delta y\,\varphi_i$ and sum over $i\in\{0,\cdots,N_y-1\}$ to get
\begin{eqnarray}
\label{eq:00}
\sum_{i=0}^{N_y-1} \Delta y\,\Q_{K,i}\,\varphi_{i} &=&  \Delta y^2\,\left( I_1 \,-\, I_2\right),
\end{eqnarray}
with 
\begin{eqnarray*}
I_1 &=& \frac{1}{2}\sum_{i=0}^{N_y-1}
\sum^{i}_{j=0} \left(a_{j,i-j}\,f_{K,j}\,f_{K,i-j}  \,-\, b_{j,i-j}\,f_{K,i}\right)\,\varphi_{i}
\\
I_2 &=& \sum_{i=0}^{N_y-1}\sum^{N_y-1}_{j=i} \left(a_{i,j-i}\,f_{K,i}\,f_{K,j-i}  \,-\, b_{i,j-i}\,f_{K,j}\right)\,\varphi_{i}.
\end{eqnarray*}
We keep the first term $I_1$ as it is and only treat the term $I_2$. First, we  commute indexes $(i,j)$ and then perform a change of index on $I_2$ which yields
\begin{eqnarray}
\nonumber
I_2 &=& \sum_{j=0}^{N_y-1}\sum^{j}_{i=0} \left(a_{i,j-i}\,f_{K,i}\,f_{K,j-i}  \,-\, b_{i,j-i}\,f_{K,j}\right)\,\varphi_{i} 
\\
\label{eq:01}
&=& \sum_{i=0}^{N_y-1}\sum^{i}_{l=0} \left(a_{i-l,l}\,f_{K,i-l}\,f_{K,l}  \,-\, b_{i-l,l}\,f_{K,i}\right)\,\varphi_{i-l}.
\end{eqnarray}
Using the symmetry  $a_{i,j}=a_{j,i}$ and $b_{i,j}=b_{j,i}$ for each $(i,j)\in\{0,\cdots,N_y-1\}$, we also have that 
\begin{equation}
\label{eq:02}
I_2 \,=\, \sum_{i=0}^{N_y-1}\sum^{i}_{l=0} \left(a_{l,i-l}\,f_{K,l}\,f_{K,i-l}  \,-\, b_{l,i-l}\,f_{K,i}\right)\,\varphi_{l}.
\end{equation}
Finally, gathering  (\ref{eq:00})-(\ref{eq:02}), we get the result
$$
\sum_{i=0}^{N_y-1} \Delta y\,\Q_{K,i}\,\varphi_{i} \,=\, -\frac{1}{2}\sum_{i=0}^{N_y-1}\sum^{i}_{l=0} \Delta y^2 \,\left( a_{l,i-l}\, f_{K,l} \, f_{K,i-l} \,-\, b_{l,i-l}\,f_{K,i} \right)\left(\varphi_{l} \,+\,\varphi_{i-l} \,-\, \varphi_{i} \right).
$$
\endproof

Now we are ready to prove Theorem~\ref{thm01}.



\section{Proof of Theorem~\ref{thm01}}
 We assume that the kinetic and diffusion coefficients fulfil (\ref{hypd}), (\ref{hypa}) and (\ref{hypM}), respectively, and that $a$ and $b$ are positive {\it a.e.} in $\RR^+ \times \RR^+$. We are also given an initial datum $f^{\texttt in}$ satisfying (\ref{hypL1}) and denote by $f^h$ the solution to (\ref{discr}) constructed in (\ref{QKi})-(\ref{discrF}).

On the one hand, existence and uniqueness of a solution  to
\begin{equation}
\label{eq:000}
 \frac{df_{K,i}}{dt}  \,-\, d(y_{i-\half})\sum_{\sigma\in \mathcal{E}_K}\tau_{\sigma} \,D_{K,\sigma}f_{K,i} \, =\,  {\rm m}(K) \,\Q_{K,i},
\end{equation}
directly follows by applying the classical Cauchy-Lipschitz theorem for a finite set of ordinary differential equations.  Existence of a global solution will be a consequence of the following uniform estimates with respect to time.

On the other hand, since the discrete distribution function is solution to (\ref{eq:000}),  the nonnegativity of $f^h(t)$ follows from the monotonicity of the discrete operator $\sum_{\sigma\in \mathcal{E}_K}\tau_{\sigma} \,D_{K,\sigma}f_{K,i}$ and the structure of the discrete coagulation-fragmentation operator (\ref{discrQ}) written as the sum of a gain operator and a local loss term.

Next, we prove the total volume conservation 
\begin{equation}
\label{esti:1}
\frac{d}{dt} \sum_{i=0}^{N_y-1} \sum_{K\in\T} {\rm m}(K)\, \Delta y \,y_{i-1/2} \,f_{K,i}(t) \,\,=\,\, 0.
\end{equation}
Indeed; we multiply (\ref{eq:000}) by $\Delta y\,y_{i-1/2}$ and sum over $(K,i)\in\T\times\{0,\ldots,N_y-1\}$
\begin{eqnarray*}
\frac{d}{dt} \sum_{K,i} {\rm m}(K)\,\Delta y \, y_{i-1/2}\,f_{K,i}(t) &=& -\frac{1}{2}\sum_{K,i}\Delta y \, d(y_{i-\half})\sum_{\sigma\in \mathcal{E}_K}\tau_{\sigma} \,D_{K,\sigma}f_{K,i} D_{K,\sigma}(y_{i-1/2})  \, 
\\
&+&  \sum_{K,i} {\rm m}(K) \,\Delta y \,\Q_{K,i} y_{i-1/2}.
\end{eqnarray*}
Thus using a discrete integration by part, we prove that the first term of the right hand side is zero, whereas applying Proposition~\ref{prop:weak} with $\varphi_i=y_{i-1/2}$, it also gives that
\begin{equation*}
\sum_{K,i} {\rm m}(K) \,\Delta y \,\Q_{K,i} y_{i-1/2}=0,
\end{equation*}
which concludes the proof of $(ii)$.

Now, let us prove the stabilization towards an equilibrium in the long time when the kinetic coefficients $a$ and $b$ satisfy the detailed balance condition (\ref{dbc}).
As mentioned in the introduction the detailed balance condition (\ref{dbc}) ensures an analogue of the Boltzmann $H$-theorem for the coagulation-fragmentation equations which we derive now for the discrete solution to (\ref{discr}). We set 
$$
H(f^h|M) \,=\, \sum_{K,i} {\rm m}(K)\, \Delta y\,\left( f_{K,i}\,\left(\ln\left(\frac{f_{K,i}}{M_i}\right)\,-\, 1  \right)  + M_i\right) 
$$
and take in the discrete weak formulation (\ref{weak:d}), the test function $\varphi$ such as $\varphi_i = \ln(f_{K,i}/M_i)$ with $M_i=M(y_{i-1/2})$ and noticing that
\begin{equation*}
\frac{dH}{dt}(f^h) \,=\, \sum_{K,i} {\rm m}(K)\, \Delta y\,\frac{d f_{K,i}}{dt}\,\ln\left(\frac{f_{K,i}}{M_i}\right), 
\end{equation*}
and recalling (\ref{weak:d}) with  $\varphi_i = \ln(f_{K,i}/M_i)$, we obtain 
\begin{eqnarray*}
&& \frac{dH}{dt}(f^h) \,+\, \Delta y\,\sum_{K,i} \,\sum_{\sigma\in \mathcal{E}_K}d(y_{i-\half}) \tau_{\sigma} \,D_{K,\sigma}f_{K,i} \ln\left(\frac{f_{K,i}}{M_i}\right)
\\
 &=&    -\Delta y\,\sum_{K,i} {\rm m}(K)\,\Q_{K,i}\, \ln\left(\frac{f_{K,i}}{M_i}\right), 
\end{eqnarray*}
which yields
\begin{equation}
\label{esti:2}
\frac{dH}{dt}(f^h)  \,+\, \Delta y\,\sum_{K,i}\,\sum_{\sigma\in \mathcal{E}_K}d(y_{i-\half})\,\tau_{\sigma} \,\frac{\left(D_{K,\sigma}f_{K,i}\right)^2}{f_{K|L,i}} \,  =\, -D(f^h),
\end{equation}
with $f_{K|L,i} = (f_{K,i}-f_{L,i})/(\ln(f_{K,i}/f_{L,i}))>0$.
Then,  we treat the right hand side using the detailed balance condition (\ref{dbc}) and get
$$
D(f^h) \,=\,  \frac{1}{2}\sum_{K,i}\sum^{i}_{l=0} \Delta y^2 \,\left( a_{l,i-l}\, f_{K,l} \, f_{K,i-l} \,-\, b_{l,i-l}\,f_{K,i} \right)\left[ \ln\left(\frac{a_{l,i-l}\, f_{K,l} \, f_{K,i-l}}{b_{l,i-l}\,f_{K,i}} \right)\right],
$$
which proves $(iii)$.

Now, we study the long time asymptotic behavior of the numerical solution $f^h$. To this aim we establish the following estimates 
\begin{lemma}
\label{lmm:01}
The numerical solution $f^h$ satisfies the following uniform estimates with respect to time $t\in\RR^+$ : there exists $\kappa_0>0$ such that 
\begin{eqnarray}
\nonumber
 \sum_{K,i} {\rm m}(K)\,\Delta y\, (1+ y_{i-1/2}) \,f_{K,i}(t) \,+\,  H(f^h(t)|M) \,\leq\, \kappa_0.
\end{eqnarray}
Moreover, for all $t\in\RR^+$, we have
\begin{equation}
\label{res:1}
\Delta y \,\int_0^t  \left[\sum_{i=0}^{N_y}\,\sum_{\sigma\in \Eint}d(y_{i-\half})^{1/2}\,\medge \,|Df_{K,i,\sigma}(s)|\,\right]^2 ds  \,\leq \, \kappa_0.
\end{equation}
and
\begin{equation}
\label{res:0}
\Delta y^2 \,\int_0^t \sum_{K,i} \sum^{i}_{l=0} \,\left( a_{l,i-l}\, f_{K,l} \, f_{K,i-l} \,-\, b_{l,i-l}\,f_{K,i} \right)\left[ \ln\left(\frac{a_{l,i-l}\, f_{K,l} \, f_{K,i-l}}{b_{l,i-l}\,f_{K,i}}\right) \right]\,ds  \,\leq\, \kappa_0,
\end{equation}
\end{lemma}
\begproof
Using the total volume conservation and the entropy dissipation, we have already proven that
$$
\sum_{K,i} {\rm m}(K)\,\Delta y\, y_{i-1/2} \,f_{K,i}(t) \,+\,  H(f^h(t)|M) \,\leq\, C_0.
$$
It remains to show that there exists a constant $C_0>0$, only depending on the initial datum $f^{\rm in}$ and $M$, such that
\begin{equation}
\label{esti:3}
\sum_{K,i} {\rm m}(K)\,\Delta y\, f_{K,i}(t) \,\leq\, C_0.
\end{equation}
Indeed, on the one hand we notice that
\begin{eqnarray*}
\sum_{K,i} {\rm m}(K)\,\Delta y\, f_{K,i}(t)  &\leq &  e^2 \,\sum_{K,i} {\rm m}(K)\,\Delta y\,\mathbf{1}_{ \{ f_{K,i}(t) \leq e^2\,M_i\} }\, M_i 
 \\
&+& \frac{1}{2} \, \sum_{K,i} {\rm m}(K)\,\Delta y\,\mathbf{1}_{ \{ f_{K,i}(t) > e^2\,M_i\} }\, f_{K,i}\left|\ln\left( \frac{f_{K,i}(t)}{M_i}\right) \right|,
\\
&\leq&   e^2 \,\sum_{K,i} {\rm m}(K)\,\Delta y\, M_i 
 \\
&+& \frac{1}{2} \, \sum_{K,i} {\rm m}(K)\,\Delta y\, f_{K,i}\left|\ln\left( \frac{f_{K,i}(t)}{M_i}\right) \right|.
\end{eqnarray*}
On the other hand, observing that $r\,\ln(r) \geq r\,|\ln(r)| - 2/e$ for $r>0$, we have
$$
f_{K,i}\,\ln\left(\frac{f_{K,i}}{M_i}\right) \,\geq \,f_{K,i}\,\left|\ln\left(\frac{f_{K,i}}{M_i}\right)\right| - \frac{2 \,M_i}{e},
$$
hence it gives
\begin{eqnarray*}
\sum_{K,i} {\rm m}(K)\,\Delta y\, f_{K,i}\left|\ln\left( \frac{f_{K,i}(t)}{M_i}\right) \right|\,&\leq& \, \sum_{K,i} {\rm m}(K)\,\Delta y\, f_{K,i}\,\ln\left( \frac{f_{K,i}(t)}{M_i}\right) 
\\
&+&  \frac{2}{e}\,\sum_{K,i} {\rm m}(K)\,\Delta y\,M_i,
\\
&\leq& \,  H(f^h|M) \,+\,\,\sum_{K,i} {\rm m}(K)\,\Delta y\,\left(f_{K,i} + \frac{2}{e}\,M_i\right).
\end{eqnarray*}
Combining the two inequalities yields
\begin{eqnarray*}
\frac{1}{2} \sum_{K,i} {\rm m}(K)\,\Delta y\,f_{K,i}(t)  &\leq & \left( e^2 + \frac{1}{e}\right)  \,\sum_{K,i} {\rm m}(K)\,\Delta y\,M_i  \,+\, \frac{1}{2}\, H(f^h(t)|M),
\\
&\leq &  \left( e^2 + \frac{1}{e}\right)  \,\sum_{K,i} {\rm m}(K)\,\Delta y\,M_i  \,+\, \frac{1}{2}\, H(f^h(0)|M) =: C_0.
\end{eqnarray*}

Next, let us prove (\ref{res:0}). We start with the entropy estimates (\ref{esti:2}), which we integrate with respect to time $t\in\RR^+$, it gives
\begin{eqnarray*}
&& H(f^h(t)|M) + \int_0^t \Delta y\,\sum_{K,i}\,\sum_{\sigma\in \mathcal{E}_K}d(y_{i-\half})\,\tau_{\sigma} \,\frac{\left(D_{K,\sigma}f_{K,i}\right)^2}{f_{K|L,i}} ds \,+\,  \int_0^t D(f^h(s))ds 
\\
&\leq & H(f^h(0)|M).
\end{eqnarray*}

Then, since $H$ is a convex function, we know that $H(f^h(t)|M)$ is bounded from below and we get the result (\ref{res:0}).
 
Finally, we show there exists $C_0>$, only depending on $f^{\rm in}$ and $M$, such that 
\begin{equation}
\label{esti:4}
\int_0^t\left[\sum_{i=0}^{N_y}\sum_{{\sigma\in\Eint}\atop{\sigma=K|L}} {\rm m}(\sigma) \, d(y_{i-1/2})^{1/2}\,|Df_{K,i,\sigma}(s)|\right]^2 \,ds  \,\leq \, C_0,\qquad \forall \, t \in \RR^+. 
\end{equation}
To this aim, we start with the Cauchy-Schwarz inequality
\begin{eqnarray*}
&&\sum_{i=0}^{N_y}\sum_{{\sigma\in\Eint}\atop{\sigma=K|L}} {\rm m}(\sigma) \, d(y_{i-1/2})^{1/2}\,|Df_{K,i,\sigma}|,
\\
&\leq & \left( \sum_{i=0}^{N_y}\sum_{{\sigma\in\Eint}\atop{\sigma=K|L}} \ts \, d(y_{i-1/2})\,\frac{|Df_{K,i,\sigma}|^2}{f_{K|L,i}} \right)^{1/2} \, \left( \sum_{i=0}^{N_y}\sum_{{\sigma\in\Eint}\atop{\sigma=K|L}} \medge\dist(x_K,x_L) \,\,f_{K|L,i} \right)^{1/2}.
\end{eqnarray*}
Applying (\ref{mesh00}) and (\ref{esti:1}), (\ref{esti:2}), it follows that
\begin{eqnarray*}
&&\left(\sum_{i=0}^{N_y}\sum_{{\sigma\in\Eint}\atop{\sigma=K|L}} \Delta y \,{\rm m}(\sigma) \, d(y_{i-1/2})^{1/2}\,|Df_{K,i,\sigma}|\right)^2
\\
&\leq & \frac{2}{\xi}\,\left( \sum_{i=0}^{N_y}\sum_{{\sigma\in\Eint}\atop{\sigma=K|L}}\Delta y \, \ts \, d(y_{i-1/2})\,\frac{|Df_{K,i,\sigma}|^2}{f_{K|L,i}} \right) \,\left( \sum_{i=0}^{N_y}\sum_{K\in\T} \Delta y \,{\rm m}(K)\,f_{K,i} \right),
\\
&\leq& \frac{2\,C_0}{\xi} \,\left( \sum_{i=0}^{N_y}\sum_{{\sigma\in\Eint}\atop{\sigma=K|L}}\Delta y \, \ts \, d(y_{i-1/2})\,\frac{|Df_{K,i,\sigma}|^2}{f_{K|L,i}} \right). 
\end{eqnarray*}
Then, for each $t\in\RR^+$ we integrate the latter inequality with respect to time over the interval $[0,t]$, and use the entropy dissipation (\ref{esti:2}) which guarantees that the right hand side is uniformly bounded with respect to time.
\endproof

Now, we prove that the numerical solution converges to a steady state. Let $(t^n)_{n\geq 0}$ be a sequence of positive real numbers such that $t^n \rightarrow +\infty$ and set $f^n(t) = f^h(t^n + t)$ for $n\geq 1$ and $t \in \RR^+$. Owing to the construction of $f^n$, it is easily seen that $f^n$ is a weak solution to (\ref{discr})-(\ref{discrboundary}) on $[0,+\infty)$ with initial datum $f^h(t^n)$. We fix $T \in \RR^+$ and infer from above that
\begin{eqnarray*}
\sup_{t\in[0,T ]} \sum_{K,i} m(K)\,\Delta y\,(1 + y_i)\, f^n_{K,i}(t)  \,+\, H(f^n(t)|M) \,\leq\, \kappa_0. 
\end{eqnarray*}
Moreover, by construction of $f^n$ and applying Lemma~\ref{lmm:01}, it implies that 
\begin{eqnarray}
\label{coucou:01}
&& \Delta y \,\int_0^T  \sum_{i=0}^{N_y}\,\sum_{\sigma\in \Eint}d(y_{i-\half})^{1/2}\,\medge \,|Df^n_{K,i,\sigma}(t)|\,dt
\\
\nonumber
&\leq&  \Delta y\, \int_{t^n}^{t^n+T}\sum_{i=0}^{N_y}\,\sum_{\sigma\in \Eint}d(y_{i-\half})^{1/2}\,\medge \,|Df_{K,i,\sigma}(t)|\,dt \,\,{{\longrightarrow}\atop{n\rightarrow \infty}} \,\, 0.
\end{eqnarray}
and
\begin{eqnarray}
\label{coucou:02}
&&\Delta y^2 \,\int_0^T \sum_{K,i} \sum^{i}_{l=0} \,\left( a_{l,i-l}\, f^n_{K,l} \, f^n_{K,i-l} \,-\, b_{l,i-l}\,f^n_{K,i} \right)\left[ \ln\left(\frac{a_{l,i-l}\, f^n_{K,l} \, f^n_{K,i-l}}{b_{l,i-l}\,f^n_{K,i}}\right) \right]\,dt 
\\
\nonumber
&\leq& \Delta y^2 \,\int_{t^n}^{t^n+T} \sum_{K,i} \sum^{i}_{l=0} \,\left( a_{l,i-l}\, f_{K,l} \, f_{K,i-l} \,-\, b_{l,i-l}\,f_{K,i} \right)\left[ \ln\left(\frac{a_{l,i-l}\, f_{K,l} \, f_{K,i-l}}{b_{l,i-l}\,f_{K,i}}\right) \right]\,dt 
\\
\nonumber
&& 
\longrightarrow\,\, 0, \quad {\rm as} \,\, n \rightarrow \infty.
\end{eqnarray}
Thanks to theses estimates, we deduce that there are a subsequence of $(f^n)_{n\in\NN}$ (not relabeled) and a weak solution $\bar{f}$ to (\ref{discr})-(\ref{discrboundary}) such that
\begin{equation}
\label{lim:fn}
f^n \rightarrow  \bar{f}, \,\, {\rm in } \, C([0, T ]),\,\,{\rm as}\, n\rightarrow +\infty.
\end{equation}
Moreover, passing to the limit in (\ref{discrQ}), it easily follows that
\begin{equation}
\Q_{K,i}(f^n) \,\longrightarrow\, \Q_{K,i}(\bar{f}), \,\, {\rm in } \,  L^1(0, T),\,\,{\rm as}\, n\rightarrow +\infty.
\label{res:00}
\end{equation}
and passing to the limit in the discrete diffusion operator, we get 
\begin{equation}
\label{res:01}
d(y_{i-\half})\sum_{\sigma\in \mathcal{E}_K}\tau_{\sigma} \,Df^n_{K,i,\sigma}\,\longrightarrow \, d(y_{i-\half})\sum_{\sigma\in \mathcal{E}_K}\tau_{\sigma} \,D\bar{f}_{K,i,\sigma} \,\,{\rm in } \, L^1(0, T ),\,\,{\rm as}\, n\rightarrow +\infty.
\end{equation}

On the one hand, we use (\ref{coucou:02}) and   (\ref{lim:fn}) to conclude
that 
$$
D(\bar{f}) \,=\, 0,\,\, a.e. \,{\rm in}\,\, (0, T ).
$$ 
Consequently, $\bar{f}$ satisfies 
$$
a_{i,j} \,\bar{f}_{K,i}(t) \,\bar{f}_{K,j}(t) \,\,=\,\, b_{i,j}\,\bar{f}_{K,i+j}(t), 
$$
for almost every $t \in (0, T )$ and each $(K,i,j)\in \T\times\{0,\ldots,N_y\}^2$. In particular, it implies that
$$
\Q_{K,i}(\bar{f}(t)) = 0, \quad (K,i,j)\in \T\times\{0,\ldots,N_y\}^2.
$$ 
On the other hand, (\ref{coucou:01}) and (\ref{res:01}) ensure that
\begin{equation}
\label{coucou}
 D_{K,\sigma}\bar{f}_{K,i} \,=\, 0, \quad (K,i,\sigma)\in \T\times\{0,\ldots,N_y\}\times\E
\end{equation}
for almost every  $t\in(0,T)$. Therefore, $\bar{f}(t)$ does not depend on time and satisfies the steady states equation (\ref{discr}). Moreover, using the zero flux boundary conditions (\ref{hypBC}) and (\ref{coucou}), $\bar{f}$ does  not depend on $K\in\T$. Recalling that the differential equation (\ref{discr}) conserves global volume, we conclude there is  $\alpha\in \RR^+$ such that
$$
\bar{f}_{K,i} \,:=\, M_i, \,\,a.e.\,\, {\rm in} \,\,(0, T ),
$$
where 
$$
\sum_{i=0}^{N_y-1} \sum_{K\in\T} \Delta y \,{\rm m}(K) \, y_{i-1/2} M_{i} \,=\,  \sum_{i=0}^{N_y-1} \sum_{K\in\T} \Delta y \,{\rm m}(K) \, y_{i-1/2} f^\texttt{in}_{K,i}.  
$$

\section{Numerical simulations}\label{test:num}
This section is devoted to the numerical study of the convergence to equilibrium under the detailed balance condition and when this condition is not satisfied. We also investigate numerically the case with non homogeneous Dirichlet boundary conditions.

\subsection{Detailed balance kernels and convergence to equilibrium}
We assume that the coagulation and fragmentation coefficients fulfil the detailed balance condition: there exists a nonnegative function $M\in L^1_1(\RR^+\times\RR^+)$, such that
\begin{equation}
\label{det:bal}
a(y,y') \,M(y)\,M(y') \,\,=\,\, b(y,y')\,M(y+y'), \quad (y,y')\,\in \,\RR^+\times\RR^+.
\end{equation}
Since in that case there exists a Lyapunov functional $H$ at the discrete level, we have proven that the numerical solution $f^h$ converges to a discrete equilibrium $M^h$. Here, we choose  kernels $a$ and $b$ as follows:
$$
a(y,y') \,=\, b(y,y') \,=\, 1,
$$
which yields $M(y) = \exp(-y/\sqrt{M_1})$. It is the Aizenman-Bak model for reacting polymers which diffuse in space with a non degenerate size-dependent coefficient $d(y)=\alpha>0$. In \cite{CDF}, the authors demonstrate that the entropy-entropy dissipation methods developed by Desvillettes and Villani for the Boltzmann equation \cite{DV} applies directly and gives the exponential convergence with explicit rates towards global equilibrium for constant diffusion coefficient in any spatial dimension or for the non degenerate diffusion in dimension one. Thuerefore, the global equilibrium is given by
$$
M(y) \,=\, \exp(-y/\sqrt{M_1}), \quad y\in\RR^+, 
$$
where the number $M_1$ is given by
$$
M_1 \,:=\, \frac{1}{|\Omega|}\,\int_\Omega\,\int_{\RR^+} y\,f^\texttt{in}(x,y)\,dx\,dy, 
$$
and $M$ satisfies 
$$
d(y)\,\Delta_x M  \,+\, \Q(M) \,=\, 0.
$$
On the other hand, the local equilibrium $M_{loc}$ is
$$
M_{loc}(x,y) = \exp(-y/\sqrt{M_1(x)})
$$
where the function $x\in\Omega \,\rightarrow\, M_1(x)\in\RR^+$ is given by
$$
M_1(x) \,:=\, \int_{\RR^+} y\,f(t,x,y)\,dy, 
$$
and $M_{loc}$ satisfies 
$$
\Q(M_{loc}) \,=\, 0.
$$
The relative entropy $H(f|M)$ can be split in two different parts
\begin{equation}
\label{HfM}
H(f|M) \,=\, H(f|M_{loc})  +  H(M_{loc}|M),  
\end{equation}
where the first term in the right-hand side represents the ``distance'' between $f$ and the local equilibrium whereas the second term evaluates the distance between the local and global equilibria and only depends on the macroscopic quantity $M_1$.
In the following we present some numerical results. As an initial datum, we take 
$$
f^\texttt{in}(x,y) \,\,=\,\, \exp(-\,\alpha(x)\,y),
$$
with $\alpha(x) = 1+0.1\cos(2\,\pi\,x_1)\,\cos(2\,\pi\,x_2)$ and $x=(x_1,x_2)\in (-1/2,1/2)^2$. The diffusion is taken to be constant $d(y)=0.1$, $R\,=\,20$, $N_y=64$ and next $128$ and finally $\Delta t=0.002$. 

Let us first mention that in \cite{FMP,FR} a similar problem {\it i.e.}; trend to equilibrium of the solution to the nonhomogeneous Boltzmann equation, is investigated numerically. It is shown that the relative entropy with respect to the local equilibrium oscillates with time when the solution becomes close to the equilibrium. Here, we will show that the situation is completely different and much simpler.
   
\begin{figure}[tb]
\begin{tabular}{cc}
\includegraphics[width=7.cm,height=7.cm]{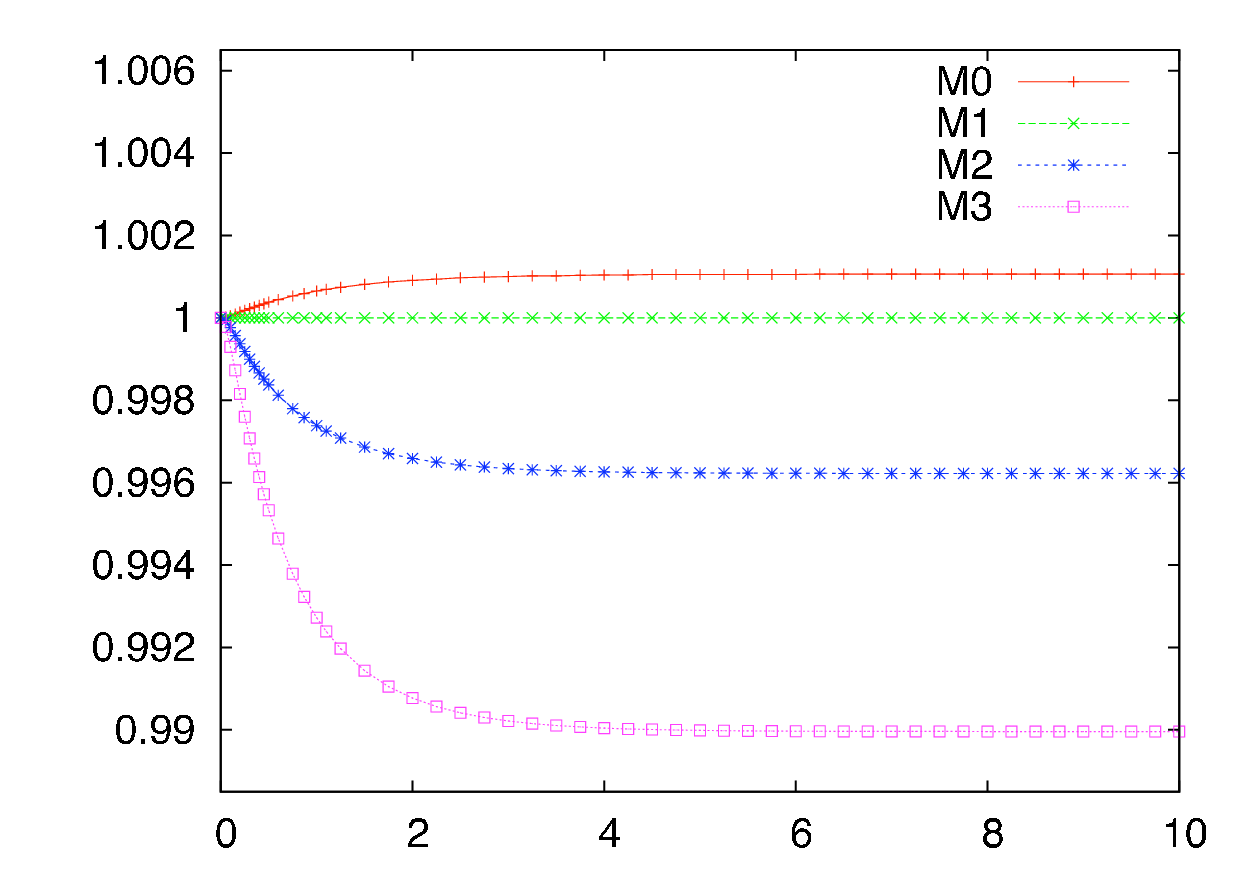}
&
\includegraphics[width=7.cm,height=7.cm]{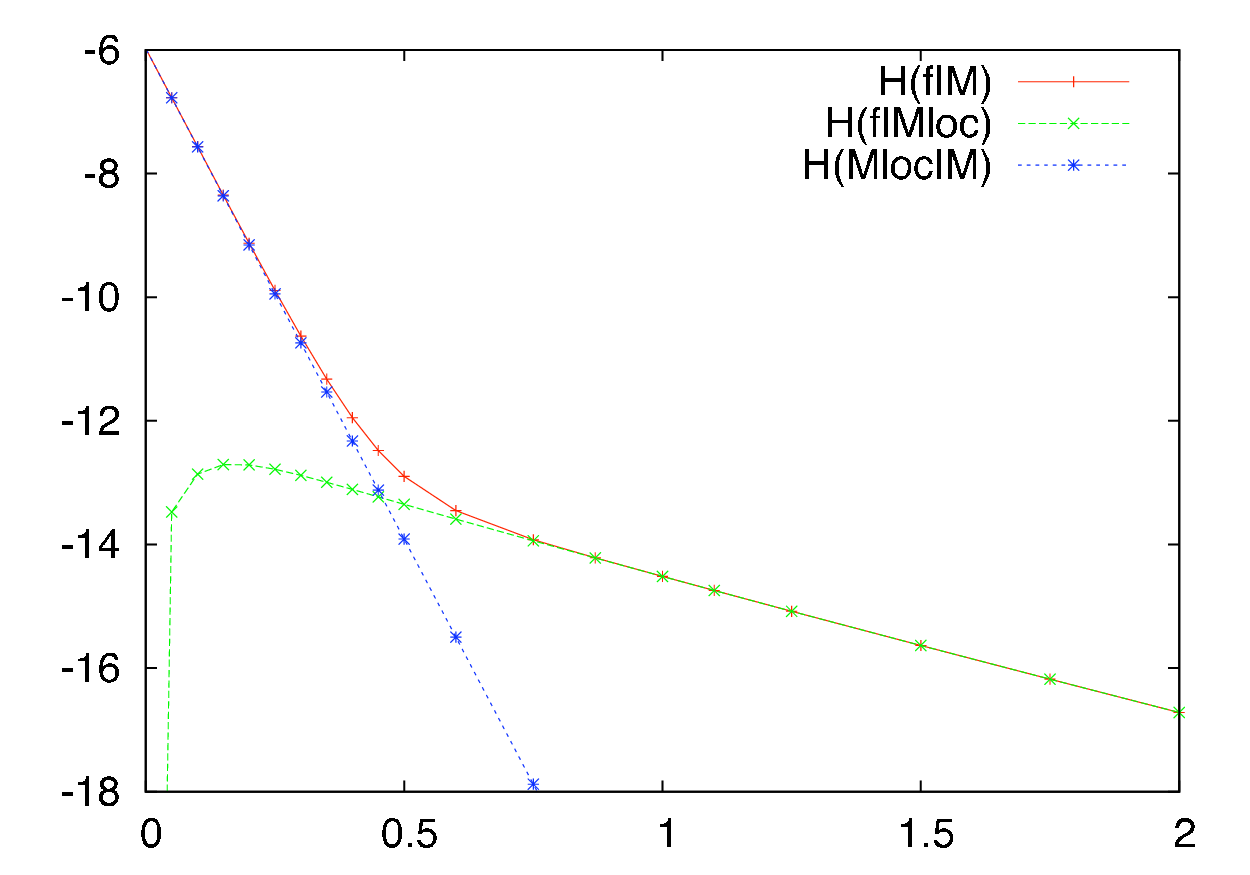}
\\
(a) & (b)
\end{tabular}
\caption{{\it Evolution of (a) the total number of particles $M_0$, the total volume $M_1$, $M_2$ and $M_3$ (b) the functional $H(f|M)$, $H(f|M_{loc})$ and $H(M_{loc}|M)$ in log scale.}}
\label{fig1}
\end{figure}

In Figure~\ref{fig1}, we report the evolution of the total number of particles $M_0$, the high order moments in $y$ of $f^h$ with respect to time and  the Lyapunov functional $H(f|M)$. As expected, the total volume $M_1(t)$ remains constant throughout time evolution and the moments stabilize to a fixed value. As regards the asymptotic profile, our numerical results are in fair agreement with the equilibrium $M(x)=\exp(-x)$. Moreover, we observe that the scheme is able to give the correct behavior of the numerical entropy $H(f|M)$, which converges exponentially fast to zero \cite{CDF}. However, we observe that trend to equilibrium for the coagulation-fragmentation with diffusion is much more simpler than trend to equilibrium for the Boltzmann equation \cite{FMP} since no oscillation occurs for the relative entropy $H(f|M_{loc})$. Indeed, after  a short transient regime, the solution $f$ converges to equilibrium as an exponential with respect to time. 

\subsection{Convergence to equilibrium for not detailed balance kernels}
\begin{figure}[tb]
\begin{tabular}{cc}
\includegraphics[width=7.cm,height=7.cm]{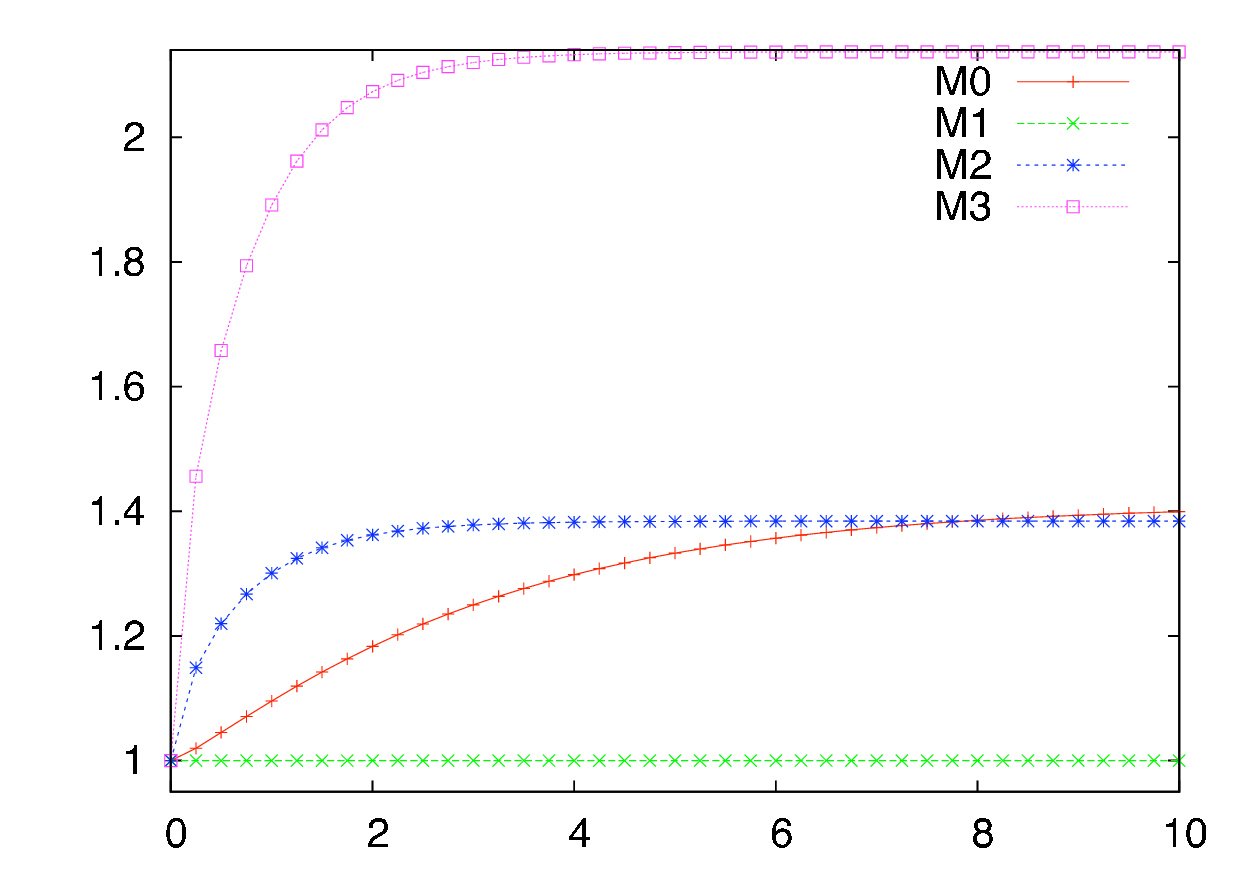}
&
\includegraphics[width=7.cm,height=7.cm]{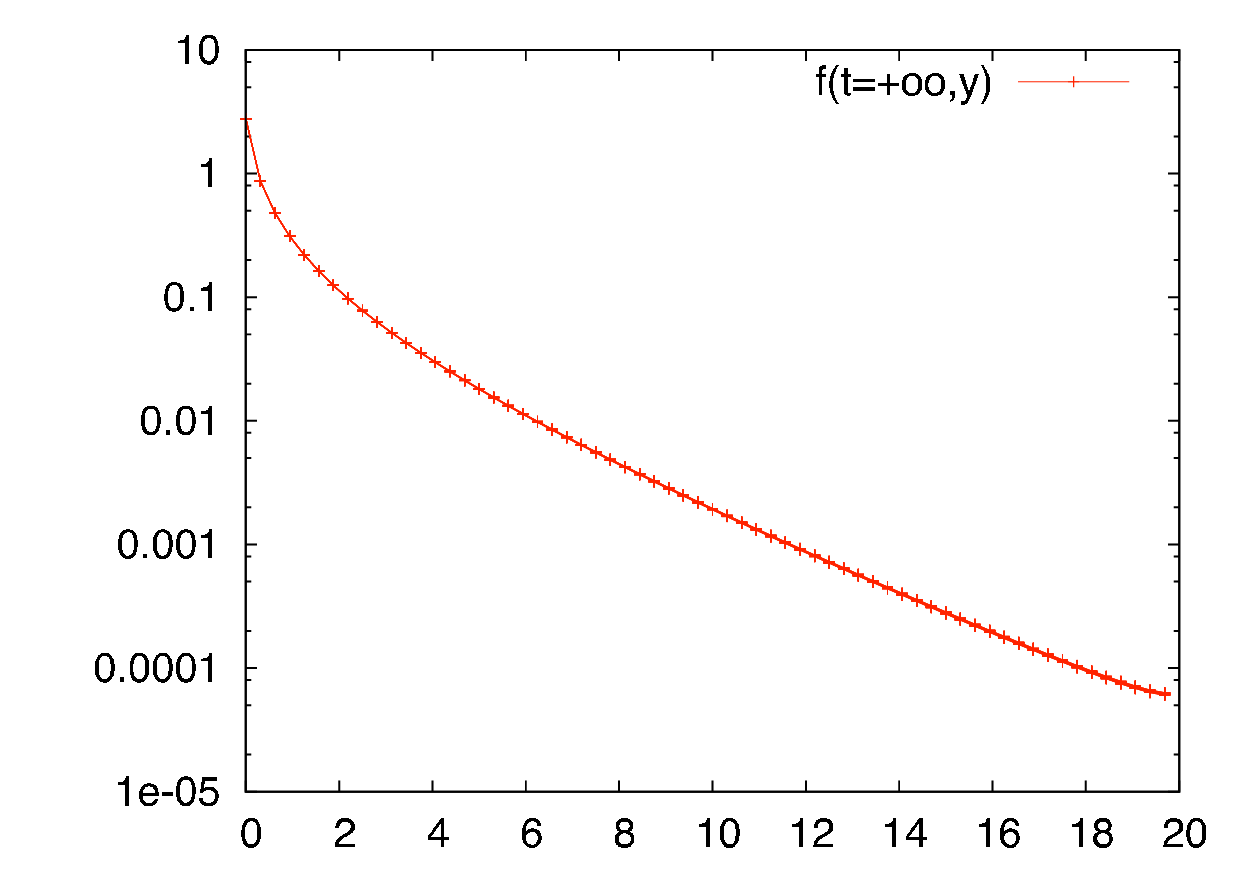}
\\
(a) & (b)
\end{tabular}
\caption{{\it Evolution of (a) the total number of particles $M_0$, the total volume $M_1$, $M_2$ and $M_3$ (b) the functional $H(f|M)$, $H(f|M_{loc})$ and $f^h(t\simeq\infty,y)$ in log scale.}}
\label{fig2}
\end{figure}

After this first result, we now explore a different situation where nothing is known about equilibrium and entropy dissipation. Indeed, we choose  kernels $a$ and $b$ as follows:
$$
a(y,y') \,=\, (y\,y')^{1/2}, \quad b(y,y') \,=\, 1.
$$
The situation becomes much more complicated since we do not know neither the expression of the steady state nor the expression of entropy. As an initial datum, we take 
$$
f^\texttt{in}(x,y) \,\,=\,\, \exp(-\,\alpha(x)\,y),
$$
with $\alpha(x) = 1+0.5\cos(4\,\pi\,x_1)\,\cos(4\,\pi\,x_2)$ and $x=(x_1,x_2)\in [-1/2,1/2]^2$. Moreover, the diffusion is taken to be 
 $$
d(y)=0.1/(1+y),
$$
which is degenerated for large $y$. We choose the truncation such that $R\,=\,20$, $N_y=64$ and next $128$ and finally $\Delta t=0.002$. In Figure~\ref{fig2}, we still report the evolution of the total number of particles $M_0$ and other high order moments in $y$ of $f^h$, that is $M_1$, $M_2$ and $M_3$.  As expected, the total volume $M_1(t)$ remains constant throughout time evolution, whereas high order moments vary and next stabilize to a fixed value, which means that the solution converges to a steady state. In the same figure (right hand side), we also report the  distribution function with respect to $y$ in log scale for large time, which indicates that the tail with respect to $y$ of the steady state $\bar{f}$ is still exponentially decreasing and $\bar{f}$ is of course constant in $x\in\Omega$.

\subsection{Convergence to equilibrium for non homogeneous boundary conditions}
In this last section, we study the coagulation-fragmentation operator with diffusion in space and mixed boundary conditions in $x\in\Omega$. More precisely, we aim to approximate the solution in $\Omega=(0,1/8)\times (0,1)$ to
$$
\left\{
\begin{array}{l}
\displaystyle{\frac{\partial f}{\partial t} \,-\, d(y)\, \Delta_x f \,= \Q(f), 
,}
\\
\,
\\
f(t=0,x,y) \,=\, f^\texttt{in}(x,y), \quad (x,y)\in \Omega \times \RR^+,
\\
\,
\\
f(t,x_1=0,x_2,y) \,=\, \exp\left(-y/\tilde{\alpha}(x_2)\right), \quad (x_2,y)\in \partial\Omega \times \RR^+,
\\
\,
\\
\nabla_x f \cdot \nu(x) \,=\, 0, \quad x_1=1/8, \,{\rm or }\,\, x_2=0, \,{\rm or} \,\, x_2=1,
\end{array}
\right.
$$
with $\tilde{\alpha}(x_2) = (1\,+\,\cos(4\pi x_2))/2$ and $\nu$ the external unit normal to $\partial\Omega$.

We consider an initial condition, which is at equilibrium for the coagulation fragmentation operator
$$
f^\texttt{in}(x,y)\,=\,exp(-y/\alpha(x)) , \quad (x,y)\in \Omega \times \RR^+,
$$
with $\alpha(x) = (1\,+\,\cos(32\,\pi\,x_1)\,\cos(4\,\pi\,x_2))/2$.

We perform numerical simulations using the finite volume scheme with diffusion coefficient $d(y)=0.01/(1+y)$ and $a=b\equiv 1$. The following figures (Fig~\ref{fig3}, ~\ref{fig4} and \ref{fig5}) are illustrations of the profile of the number of particles $M_0(t,x)$ and total volume $M_1(t,x)$  for $t\in \RR^+$ and $x\in \Omega$ and also the projections
$$
P(t,x_2,y) = \int_0^{1/8} y\,f(t,x,y)dx_1, \quad x_2\in (0,1), y\in\RR^+.
$$ 
It allows us to observe the qualitative behavior of the system. The profiles are in that case smooth and also converges to an equilibrium. For such a computations we have considered $N_y=64$ and $128\times 128$ grid points for the space variable $x\in\Omega$.  In this situation, we cannot characterize explicitly the equilibrium in space since the diffusion coefficient is not constant and $H$ is not valid here due to the Dirichlet boundary conditions for $x_1=0$ and then the total volume is not preserved at all. However, we still observe that the numerical solution converges to a discrete equilibrium (see Fig. \ref{fig3}-\ref{fig5}).   

\begin{figure}[htbp]
\begin{tabular}{cc}
\includegraphics[width=7.cm,height=7.cm]{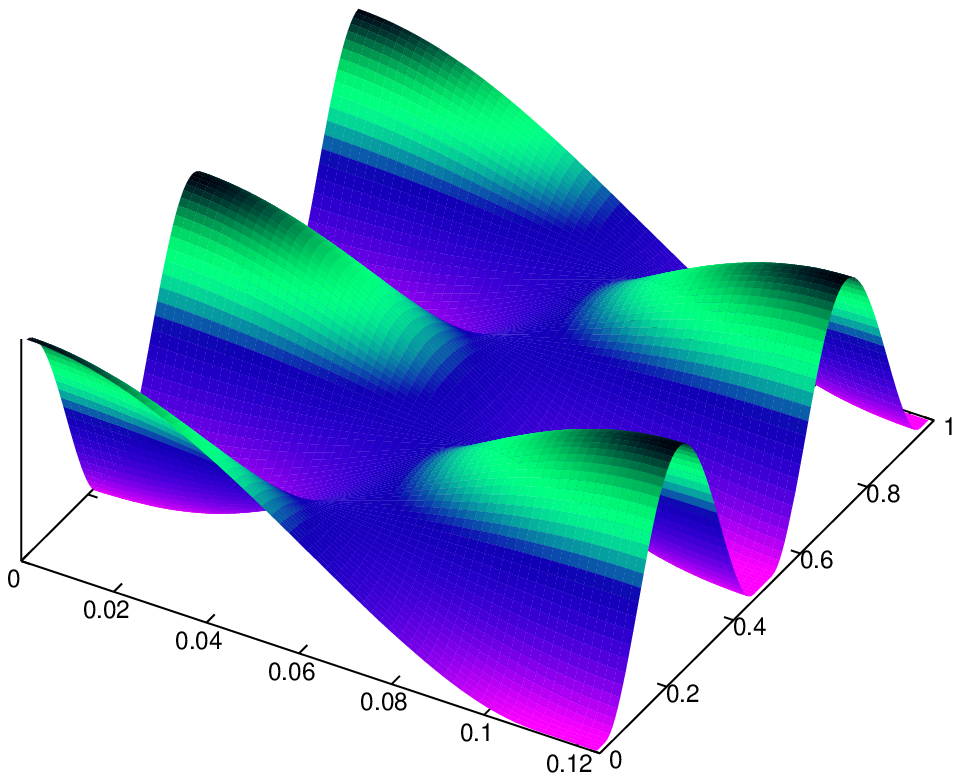}&
\includegraphics[width=7.cm,height=7.cm]{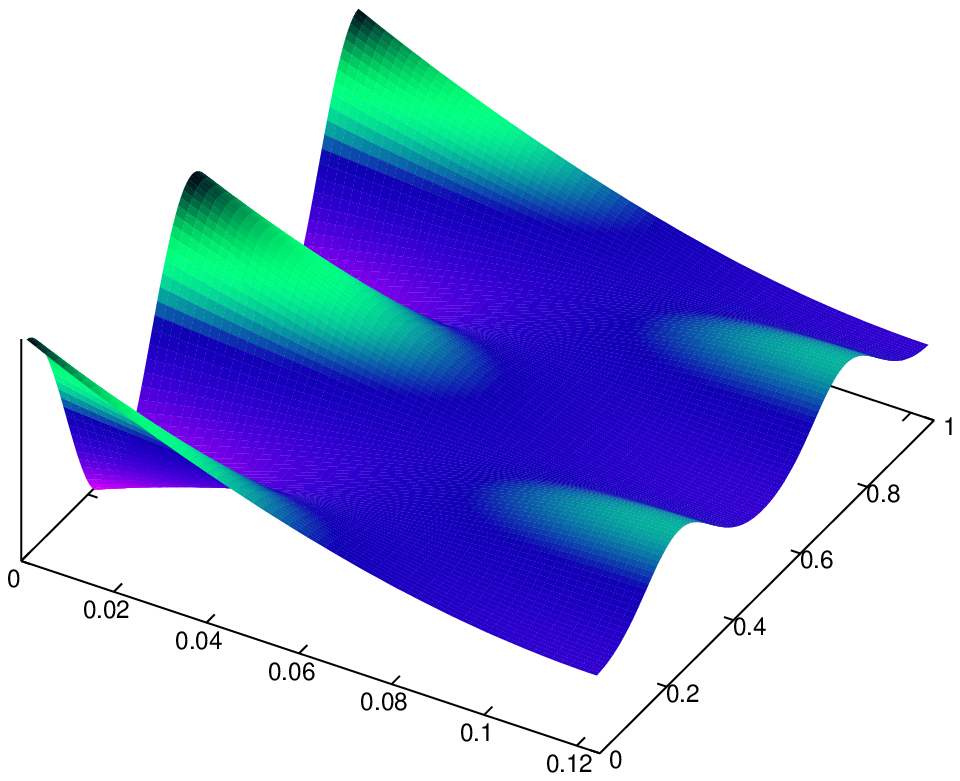}
\\
\includegraphics[width=7.cm,height=7.cm]{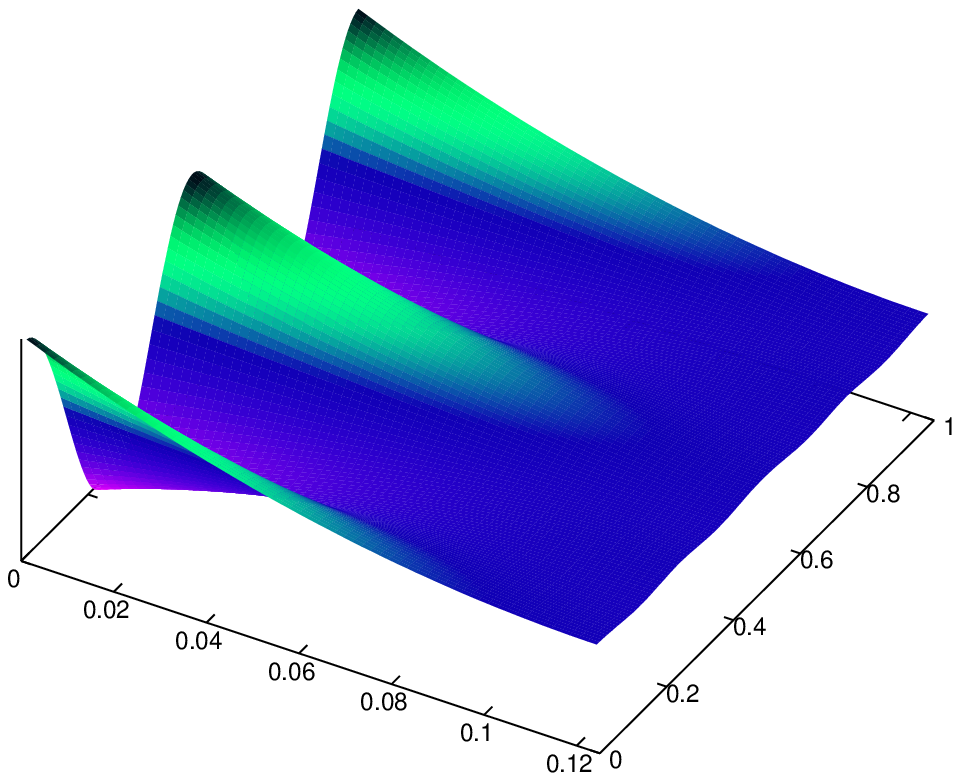}&
\includegraphics[width=7.cm,height=7.cm]{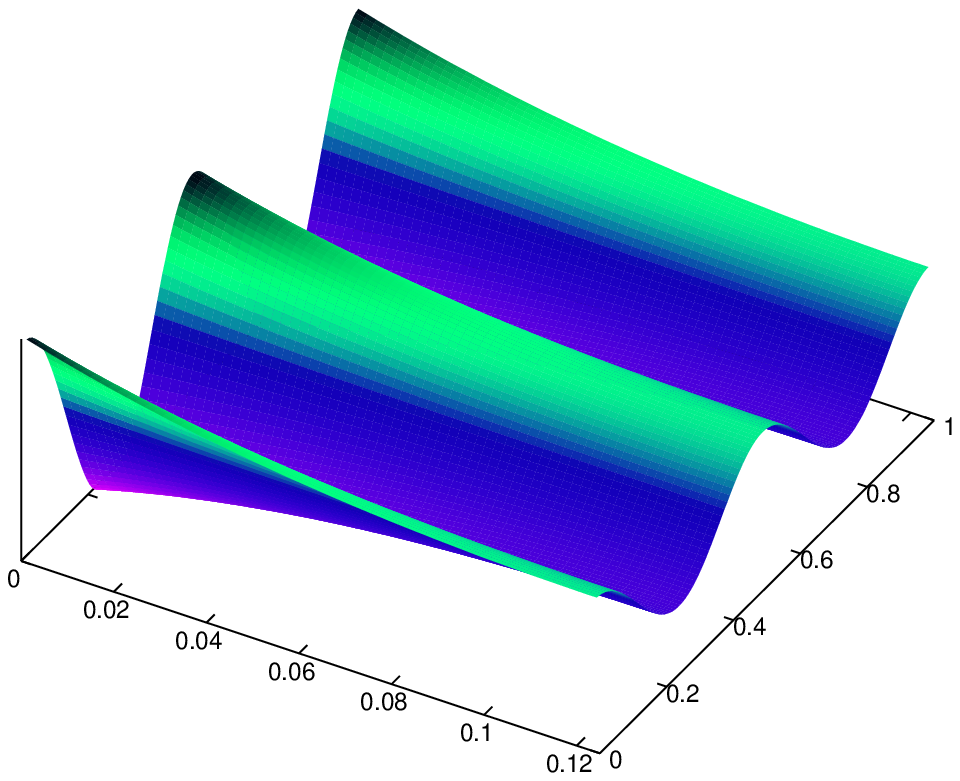}
\end{tabular}
\caption{Time evolution of the density $M_0(x)$ at time $t=0$; $0.33$; $0.66$ and $4$}
\label{fig3}
\end{figure}

\begin{figure}[htbp]
\begin{tabular}{cc}
\includegraphics[width=7.cm,height=7.cm]{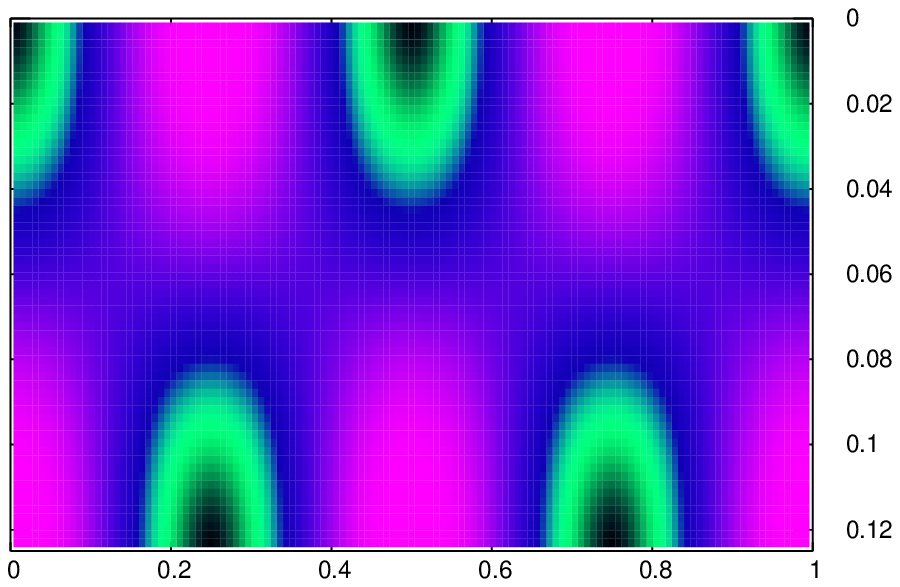}&
\includegraphics[width=7.cm,height=7.cm]{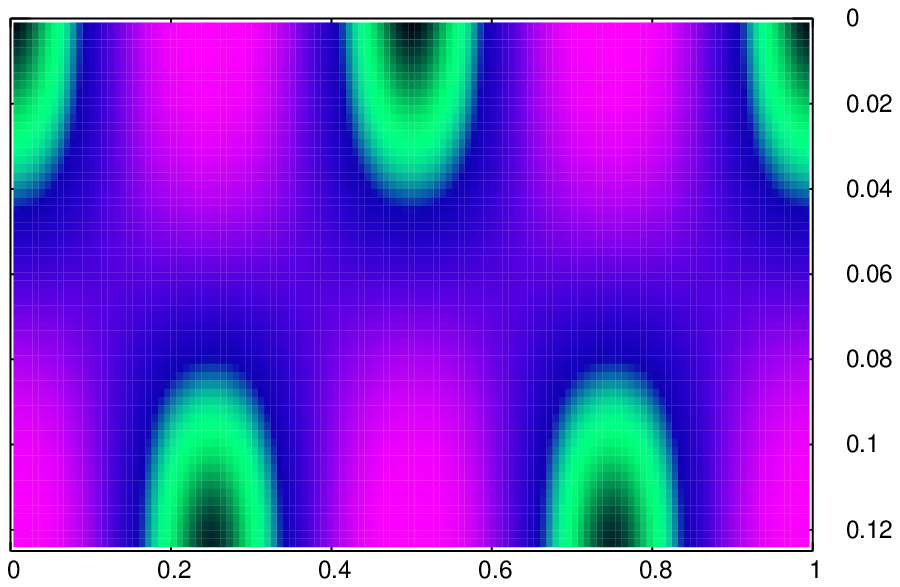}
\\
\includegraphics[width=7.cm,height=7.cm]{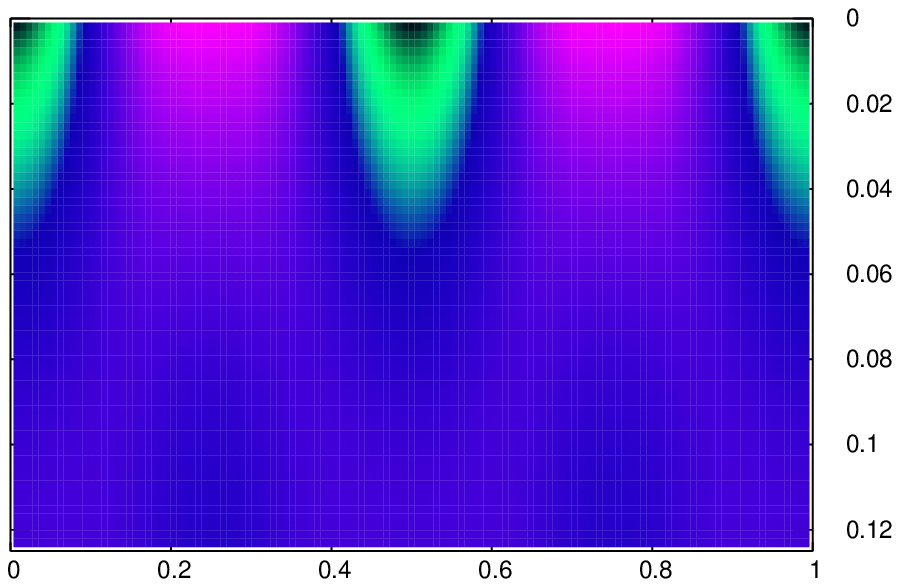}&
\includegraphics[width=7.cm,height=7.cm]{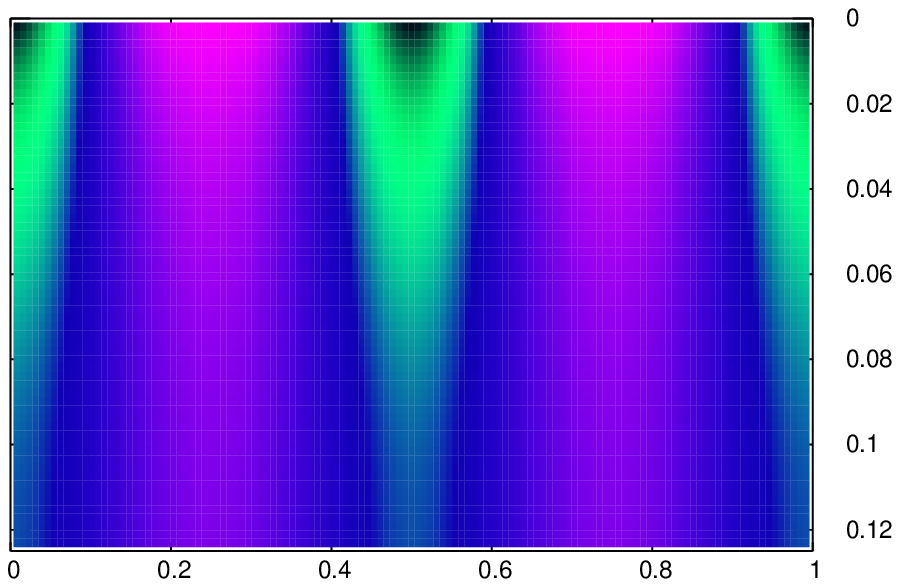}
\end{tabular}
\caption{Time evolution of the volume density $M_1(x)$ at time $t=0$; $0.33$; $0.66$ and $4$}
\label{fig4}
\end{figure}

\begin{figure}[htbp]
\begin{tabular}{ccc}
\includegraphics[width=7.cm,height=7.cm]{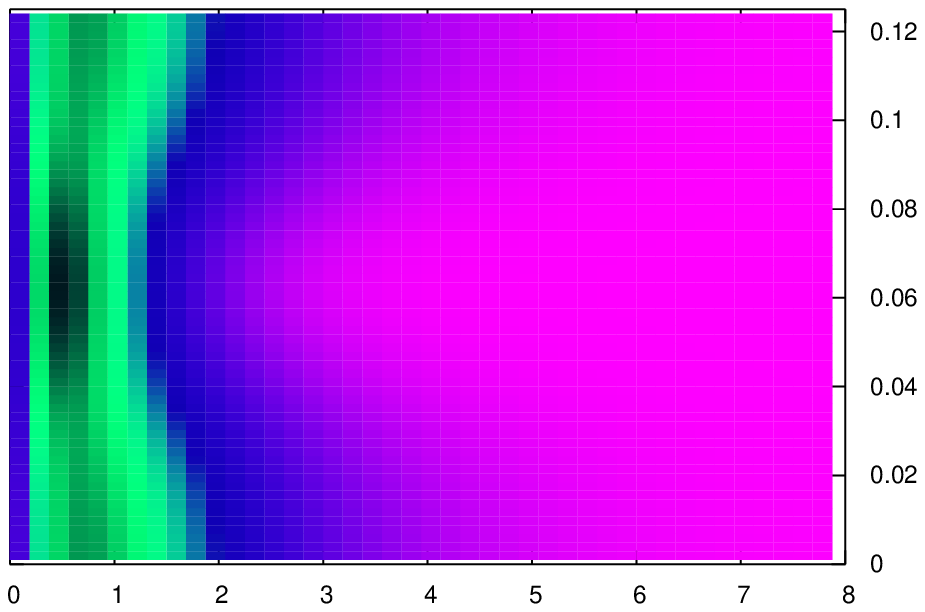}&
\includegraphics[width=7.cm,height=7.cm]{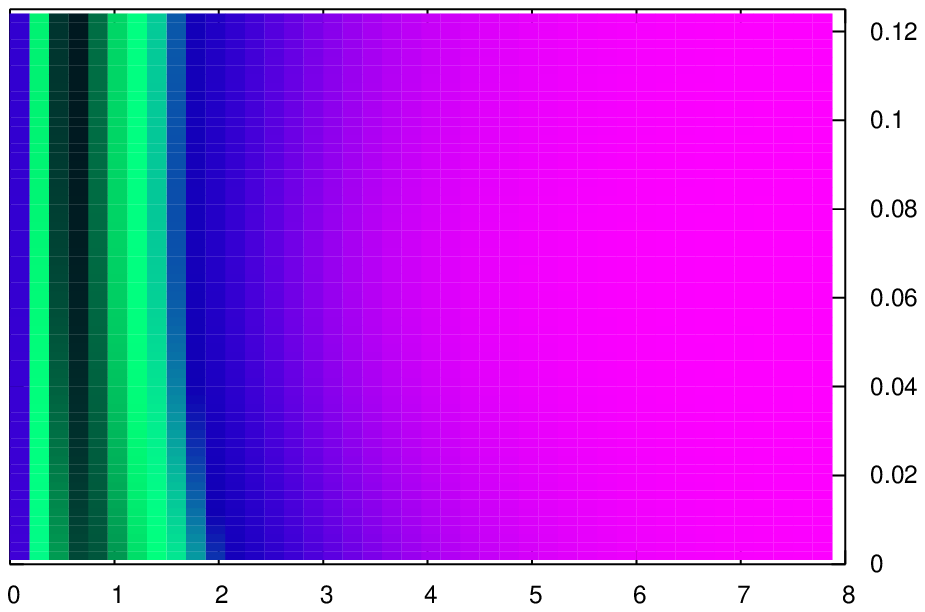}
\\
\includegraphics[width=7.cm,height=7.cm]{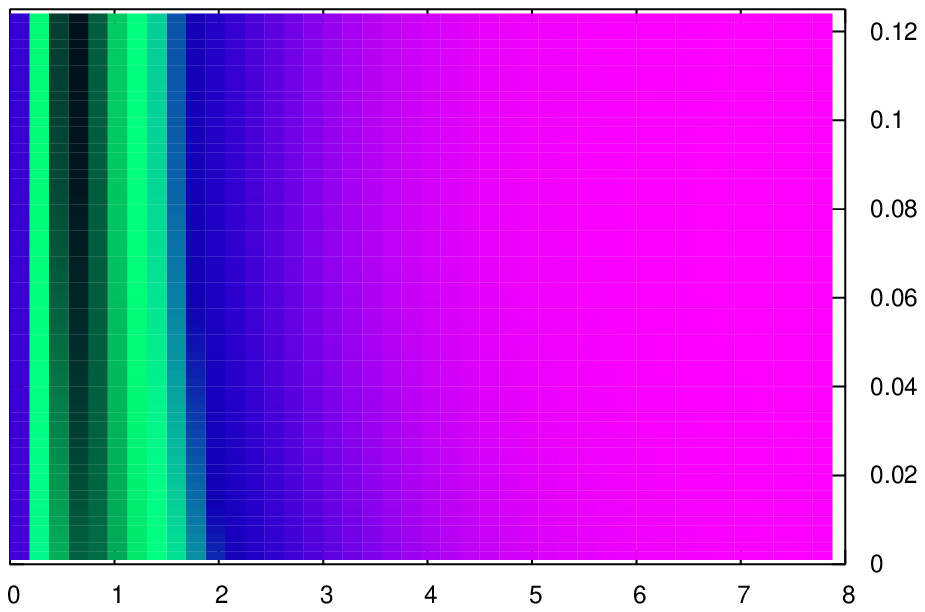}&
\includegraphics[width=7.cm,height=7.cm]{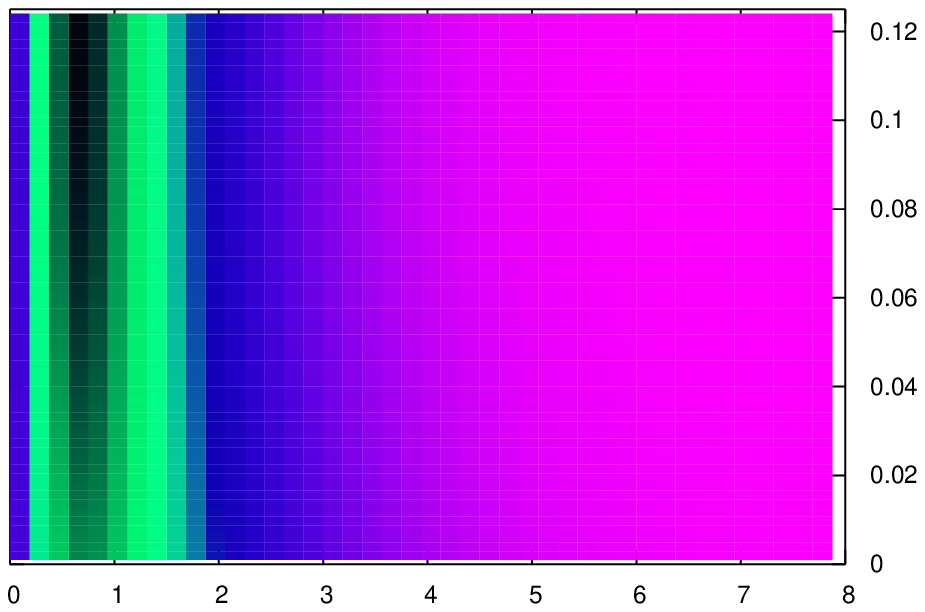}
\end{tabular}
\caption{2D projection $(x_1,y)$  of the density $y \,f(t,x,y)$ at time $t=0$; $0.33$; $0.66$ and $4$}
\label{fig5}
\end{figure}


\section{Conclusion}
In this paper, we make use of different principles (volume conservation, entropy dissipation, existence of steady states) which allow to build a stable and accurate numerical scheme for the nonlinear dynamics of coagulation, fragmentation and diffusion equations. Such an algorithm  is able to recover the main properties of the exact solution and in some particular cases, it is proven that the method is asymptotically stable, in the sense that the numerical solution converges to an equilibrium which is consistent with the exact equilibrium of the system. Numerical simulations illustrate the efficiency of the algorithm even when we cannot prove convergence to equilibrium.

 Here, we have only considered the time evolution of a distribution function, but the method can be easily coupled with Euler or Navier-Stokes equations after some adaptations. Typical applications are transport problems (including linear and nonlinear diffusion or fluid dynamics), which are coupled with population balance dynamics represented by a distribution function $f$ depending on space $x\in \Omega$ and  ``size'' variable $y\in\RR^+$. 

On the other hand, for some coagulation and fragmentation kernels \cite{FLxx,bourg:fil}, total volume is not conserved at all and then a non conservative formulation can be used to dicretize the coagulation and fragmentation operators using the same kind of finite volume scheme. Then, the finite volume approach allows to use non uniform meshes for the volume variable $y\in \RR^+$, which is particularly well suited in this case \cite{FLxx}.

\vskip 1.cm
\paragraph{Acknowledgement.}  This work is partially supported by the programm ``ANR Jeunes Chercheurs''  JCJC-0136 (projet {\it MNEC}) from the french Ministery of Research.

\begin{flushright} 
\signff 
\end{flushright}

\end{document}